\documentclass[11pt,letterpaper]{article}
\usepackage{MathStyle}
\usepackage[compact]{PaperStyle} 

\DeclareMathOperator{\disc}{disc}               
\DeclareMathOperator{\herdisc}{herdisc}         
\DeclareMathOperator{\sech}{sech}               
\DeclareMathOperator{\depth}{depth}             
\DeclareMathOperator{\lev}{lev}                 
\DeclareMathOperator{\range}{range}             
\DeclareMathOperator{\shape}{shape}             

\newcommand{\Ad}{\mathcal{A}_d}                 
\newcommand{\Bd}{\mathcal{B}_d}                 
\newcommand{\Jst}{J^{*}}                        
\newcommand{\Plo}{P_{\mathrm{lo}}}              
\newcommand{\Phig}{P_{\mathrm{hi}}}             
\newcommand{\Wlin}{W_{\mathrm{lin}}}            
\newcommand{\Flin}{F_{\mathrm{lin}}}            
\newcommand{\ent}{\mathsf{h}}                   
\newcommand{\gm}{\gamma_M}                      
\newcommand{\Shapes}{\mathcal{S}_h}             
\newcommand{\Shell}{\mathcal{R}}                
\newcommand{\Bl}[1]{\mathcal{C}_{#1}}           
\newcommand{\Low}{\mathcal{L}}                  
\newcommand{\VH}{\mathcal{V}_H}                 
\newcommand{\Gram}{\mathcal{G}}                 
\newcommand{\dy}{\mathcal{D}}                   
\newcommand{\Rad}{\varrho}                      
\newcommand{\cube}{[0,1]^d}                     
\newcommand{\Hlo}{\mathcal{H}_{\mathrm{lo}}}    
\newcommand{\cellmon}{\mathcal{M}}              

\title{Resolving Tusn\'ady's Problem}
\author{Nicola Wengiel}
\affiliation{Department of Mathematics, FAU Erlangen-N\"urnberg}
\email{nicola.wengiel@fau.de}
\hypersetup{pdftitle={Resolving Tusnady's Problem},pdfauthor={Nicola Wengiel}}

\begin{document}
\maketitle
\begin{abstract}
We prove a lower bound of $\Omega_d((\log N)^{d-1/2})$ for the worst-case combinatorial discrepancy of $N$-point sets with respect to anchored boxes in every fixed dimension $d\ge3$. Together with Nikolov's upper bound of the same order and Wei's planar lower bound, this settles Tusn\'ady's problem: the discrepancy of axis-parallel boxes is $\Theta_d((\log N)^{d-1/2})$ in every fixed dimension $d\ge2$. The proof combines smoothed product test functions, sampling discretization, and an entropy argument for conditioned Gaussian vectors. The matching upper and lower bounds are also formalized in Lean~4.
\end{abstract}
\makepapercontents

\section{Introduction}
\label{sec:Introduction}

Let $P\subset\R^d$ be a finite set and let $\calF$ be a family of subsets of $\R^d$. The \emph{combinatorial discrepancy} of $P$ with respect to $\calF$ is
\begin{equation*}
    \disc(P,\calF)\coloneqq\min_{\chi\colon P\to\{-1,1\}}\ \sup_{F\in\calF}\ \Bigl|\sum_{p\in P\cap F}\chi(p)\Bigr|.
\end{equation*}
We write $\disc(N,\calF)$ for its supremum over all $N$-point sets $P\subset\R^d$. Two families are of interest here: the family $\Bd$ of all axis-parallel boxes $\prod_{j=1}^d[a_j,b_j]$, and the family $\Ad$ of \emph{anchored boxes} $A(x)\coloneqq\{y\in\R^d:y\le x\}$, $x\in\R^d$, where $y\le x$ is understood coordinatewise. Every set $P\cap A(x)$ is of the form $P\cap B$ for a box $B\in\Bd$, and conversely, since $P$ is finite, the indicator of every set $P\cap B$ with $B\in\Bd$ is a signed sum of indicators of $2^d$ sets $P\cap A(x)$ by inclusion--exclusion, so that
\begin{equation}\label{eq:AnchoredVsBoxes}
    \disc(N,\Ad)\le\disc(N,\Bd)\le2^d\disc(N,\Ad).
\end{equation}
Discrepancy bounds for axis-parallel boxes have applications to numerical integration, orthogonal range searching, and differentially private range counting; see \cite[Section~2]{Nikolov17}.

\paragraph{Tusn\'ady's problem.}
In the early 1980s, G\'abor Tusn\'ady asked for the order of magnitude of $\disc(N,\mathcal B_2)$, and in particular whether it is bounded by a constant independent of $N$; see \cite{Beck81}, \cite[p.~125]{Matousek99book} and \cite[Section~1.7]{BansalNikolov26}. Beck \cite{Beck81} answered the latter question negatively by transferring classical lower bounds for the unnormalized Lebesgue-measure discrepancy to combinatorial discrepancy for infinitely many $N$ \cite[Proposition~1.8]{Matousek99book}; for an $N$-point set in the unit cube, this is the largest deviation, over all boxes, of the number of points in a box from $N$ times its volume. The minimum of this discrepancy over $N$-point sets in the unit cube is of order $\log N$ in the plane, the lower bound being Schmidt's theorem \cite{Schmidt72}, and at least of order $(\log N)^{(d-1)/2}$ in dimension $d$ by Roth's theorem \cite[p.~74, footnote]{Roth54}. Determining the order of $\disc(N,\Bd)$ for fixed $d\ge2$ has become known as \emph{Tusn\'ady's problem}.

Beck also proved the upper bound $O((\log N)^4)$ in the plane; this bound and its higher-dimensional counterparts were improved in a sequence of works \cite{Beck89,Bohus90,Srinivasan97,Matousek99,Larsen14,BansalGarg17}, and the currently best upper bound, due to Nikolov \cite[Theorem~1]{Nikolov17}, is
\begin{equation}\label{eq:NikolovBound}
    \disc(N,\Ad)=O_d\bigl((\log N)^{d-1/2}\bigr)\qquad\text{for every fixed }d\ge2.
\end{equation}
On the lower bound side, Matou\v{s}ek and Nikolov \cite{MatousekNikolov15} and Matou\v{s}ek, Nikolov and Talwar \cite{MNT20} used the $\gamma_2$ factorization norm to prove $\disc(N,\Ad)=\Omega_d((\log N)^{d-1})$. In the plane, Wei's preprint \cite[Theorem~1.1]{Wei26} gives a lower bound of order $(\log N)^{3/2}$, matching \eqref{eq:NikolovBound} for $d=2$. Thus, before the present work, the order of $\disc(N,\Bd)$ was known only up to a factor of order $\sqrt{\log N}$ in every fixed dimension $d\ge3$.

\subsection{Main result and consequences}
\label{sec:MainResult}

We close this gap in every fixed dimension $d\ge3$.

\begin{theorem}
    \label{thm:Main}
    For every integer $d\ge3$ there are constants $c_d>0$ and $N_d\in\N$ such that for every integer $N\ge N_d$,
    \begin{equation*}
        \disc(N,\Ad)\ge c_d(\log N)^{d-1/2}.
    \end{equation*}
\end{theorem}

\begin{corollary}
    \label{cor:Main}
    For every fixed $d\ge3$, $\disc(N,\Ad)$ and $\disc(N,\Bd)$ are both $\Theta_d\bigl((\log N)^{d-1/2}\bigr)$.
\end{corollary}

\begin{proof}
    Combine \cref{thm:Main} with \eqref{eq:NikolovBound} and \eqref{eq:AnchoredVsBoxes}.
\end{proof}

In concurrent work, Lei \cite[Theorem~1.1]{Lei26} obtains the same sharp order in every fixed dimension $d\ge2$ and shows that, with high probability, the lower bound for a uniform random sample survives the deletion of any set containing at most a sufficiently small fixed fraction of the sample.

By \cref{cor:Main}, the lower bound $\Omega_d((\log N)^{d-1})$ of \cite{MatousekNikolov15,MNT20} is not sharp, although $(\log N)^{d-1}$ had long been conjectured to be the correct order; see Srinivasan \cite[p.~694]{Srinivasan97}, who refers to Beck and S\'os \cite{BeckSos95}, and also \cite[Conjecture~1.4.1]{Bansal22} and \cite[Section~2.2]{AistleitnerBilykNikolov18}. Together with Wei's planar lower bound \cite{Wei26}, \cref{cor:Main} settles Tusn\'ady's problem in every fixed dimension $d\ge2$. Since there are $N$-point sets in $\cube$ whose Lebesgue-measure discrepancy with respect to anchored boxes is $O_d((\log N)^{d-1})$ \cite[Theorem~2.4]{Matousek99book}, \cref{thm:Main} also shows that for $d\ge3$ the largest possible combinatorial discrepancy of $N$-point sets with respect to boxes exceeds the smallest possible Lebesgue-measure discrepancy of $N$-point sets by a factor of order at least $\sqrt{\log N}$, so that Beck's transference principle is not sharp in any fixed dimension $d\ge3$.

We treat the case $d\ge3$; for these dimensions, the relative linearization error is $O_d((\log N)^{-(d-2)})$ and tends to zero as $N\to\infty$. The constants $c_d$ and $N_d$ depend on $d$; we make no attempt to optimize them.

\Cref{thm:Main} extends from boxes to the homothetic copies $\mathcal T_B\coloneqq\{tB+x:t>0,\ x\in\R^d\}$ of a fixed convex polytope $B\subset\R^d$: the lower bound follows from it by a transfer reduction, which we carry out in \cref{sec:Polytopes}, and the upper bound is Nikolov's \cite[Theorem~2]{Nikolov17}.

\begin{corollary}[Fixed convex polytopes]
    \label{cor:Polytopes}
    For every fixed $d\ge3$ and every fixed convex polytope $B\subset\R^d$ with nonempty interior,
    \begin{equation*}
        \disc(N,\mathcal T_B)=\Theta_{d,B}\bigl((\log N)^{d-1/2}\bigr).
    \end{equation*}
\end{corollary}

The lower-bound constants depend only on $d$. Together with the lower bound of \cite{MatousekNikolov15,MNT20}, the transfer reduction already confirms both conjectures stated after \cite[Theorem~3]{Nikolov17} in every dimension $d\ge2$, and \cref{cor:Polytopes} gives the stronger exponent $d-\tfrac12$ for $d\ge3$; see \cref{sec:Polytopes}.

The upper bound in \cref{cor:Main} is also attained algorithmically: a randomized offline algorithm with polynomial running time colors any $N$-point set with discrepancy $O_d((\log N)^{d-1/2})$ for anchored boxes, with high probability (\cref{prop:Algorithm}). This follows from a recent prefix-balancing theorem of Aden-Ali \cite[Theorem~1]{AdenAli26} via the reduction in Nikolov's proof of \cite[Theorem~1]{Nikolov17}, which we write out in \cref{sec:Algorithm}; by \cref{thm:Main}, the resulting discrepancy guarantee is optimal up to the constant for $d\ge3$. The failure probability, the number of operations and the model of computation are specified with the proof.

Our test-function estimates also apply to small deviations of the Brownian sheet.

\begin{corollary}[Brownian-sheet small deviations]
    \label{cor:BrownianSheet}
    Let $B_d$ be the standard Brownian sheet on $\cube$. For every fixed integer $d\ge2$,
    \begin{equation*}
        -\log\Prob(\norm{B_d}_\infty\le\delta)\asymp\delta^{-2}\bigl(\log(1/\delta)\bigr)^{2d-1},\qquad\delta\downarrow0.
    \end{equation*}
\end{corollary}

For $d\ge3$, we prove the lower bound for $-\log\Prob(\norm{B_d}_\infty\le\delta)$ in \cref{sec:BrownianSheet}; an earlier proof via metric entropy is given in \cite[Theorem~4.4]{NeumayerUllrich26} (see also \cite[Theorem~1.1]{BilykNeumayerUllrich26}). The upper bound is due to Dunker, Linde, K\"uhn and Lifshits \cite{DunkerLindeKuhnLifshits99}, and the case $d=2$ to Talagrand \cite{Talagrand94}.

\paragraph{Formal verification.}
A Lean~4 formalization proves \cref{thm:Main}, Nikolov's upper bound \eqref{eq:NikolovBound}, and the comparison \eqref{eq:AnchoredVsBoxes} between anchored and general boxes. For the planar lower bound it reuses Wei's Lean formalization \cite{WeiLean26}. Together these results extend \cref{cor:Main} to every fixed dimension $d\ge2$. The combined theorem \texttt{Tusnady.Integrated.main\char`\_sharp\char`\_order} has no hypotheses beyond $d\ge2$ and depends only on the standard axioms \texttt{propext}, \texttt{Classical.choice}, and \texttt{Quot.sound}. The source code and verification instructions are archived on Zenodo \cite{WengielLean26}. The extensions and further consequences in \cref{sec:Consequences} are outside the scope of the formalization.

\paragraph{Notation.}
From \cref{sec:Preliminaries} to \cref{sec:Conclusion}, $d\ge3$ is fixed; in \cref{sec:Consequences} the range of $d$ is stated with each result. Constants $c_d,C_d>0$ depend only on $d$ and may change from occurrence to occurrence, except where a statement fixes them (as $c_d$ and $N_d$ in \cref{thm:Main}); they never depend on the depth $h$. We write $X\lesssim Y$ or $X=O_d(Y)$ if $X\le C_dY$, $X=\Omega_d(Y)$ if $Y\lesssim X$, and $X\asymp Y$ or $X=\Theta_d(Y)$ if both hold. For $x\in\R^d$, $|x|_1\coloneqq\sum_j|x_j|$; for $p\in\cube$, $[p,1]^d\coloneqq\prod_{j=1}^d[p_j,1]$; $[n]\coloneqq\{1,\dots,n\}$; and $\inp fg\coloneqq\int fg$, over the cube or over $[0,1]$. Logarithms are natural, $\N=\{1,2,\dots\}$, $\N_0=\N\cup\{0\}$, $3\N=\{3,6,\dots\}$ and $3\N_0=3\N\cup\{0\}$. The standard Gaussian measure on $\R^M$ is $\gm$, and $\norm{\cdot}_2$ denotes both the Euclidean norm and the norm of $L^2(\cube)$. For symmetric matrices, $X\preceq Y$ means that $Y-X$ is positive semidefinite, and $\norm{\cdot}_{\mathrm{op}}$ is the operator norm with respect to the Euclidean norms.

\subsection{Overview of the proof}
\label{sec:Overview}

We fix $d\ge3$ and describe the mechanism of the proof, suppressing all constants that depend only on $d$. The construction is parametrized by a depth $h\in3\N$; the point set has $n\asymp(h+1)^{d-1}2^h$ points, so that $h\asymp\log n$, and we show that its \emph{hereditary} discrepancy for anchored boxes, that is, the largest discrepancy of one of its subsets, is at least $c_d(h+1)^{d-1/2}$ (\cref{prop:HerdiscLowerBound}). Hence some subset has discrepancy at least $c_d(h+1)^{d-1/2}$, and padding this subset with points outside all relevant anchored boxes yields \cref{thm:Main} for every large $N$.

\paragraph{Duality.}
For points $P_1,\dots,P_n\in(0,1)^d$ and a weight vector $\eta\in\R^n$, let $G_\eta(x)\coloneqq\sum_{i\le n}\eta_i\one_{\{P_i\le x\}}$, $x\in\cube$, be the weighted count of the anchored box $A(x)$. A coloring is a weight vector $\eta\in\{-1,1\}^n$, and its discrepancy is $\norm{G_\eta}_\infty\coloneqq\sup_{x\in\cube}|G_\eta(x)|$. For an integrable function $W$ on the cube, Fubini's theorem and H\"older's inequality give
\begin{equation}\label{eq:DualityIntro}
    \inp{G_\eta}{W}=\sum_{i=1}^n\eta_i\,(\Jst W)(P_i)\le\norm{G_\eta}_\infty\norm{W}_1,
    \qquad
    (\Jst W)(p)\coloneqq\int_{[p,1]^d}W(x)\,dx .
\end{equation}
A test function $W$ with $\norm{W}_1\le1$ and a large value of $\sum_i\eta_i(\Jst W)(P_i)$ thus certifies a large weighted box count. The difficulty is that the coloring is chosen adversarially after the point set, so that the test function must be adapted to it; we obtain this adaptation from a probabilistic reformulation in which the adversarial coloring is replaced by a Gaussian weight vector conditioned on a small-ball event. Concretely, we construct a family of test functions $W(u)$, $u\in[-1,1]^M$, and evaluate it at $u=U(Z)$, the Rosenblatt coordinates of the conditioned Gaussian vector $Z$: if the hereditary discrepancy is small, the conditioning event has probability at least $\mathrm{e}^{-O(M)}$ (Step 1); this forces $Z$ to correlate with $U(Z)$ (Step 2); and the test functions are designed so that this correlation produces a large weighted box count on average, contradicting the conditioning (Steps 3--5).

\paragraph{Step 1: hereditary discrepancy and Gaussian small balls.}
Let $D$ be the hereditary discrepancy of the point set for anchored boxes, and let $\eta$ be a centered Gaussian vector in $\R^n$ with all coordinate variances at most one. \Cref{lem:GaussianCover} shows that
\begin{equation}\label{eq:CoverIntro}
    \Prob(\norm{G_\eta}_\infty\le2D)\ge(p_0/2)^n,\qquad p_0\coloneqq\Prob(|g|\le1),
\end{equation}
where $g$ is a standard Gaussian variable, by combining \v{S}id\'ak's inequality \cite{Sidak67}, the rounding lemma of Lov\'asz, Spencer and Vesztergombi \cite{LSV86}, and Anderson's inequality \cite{Anderson55}.

\paragraph{Step 2: conditioned Gaussian vectors and Rosenblatt coordinates.}
We take $\eta=Bz$, where $z$ is a standard Gaussian vector in $\R^M$, $M$ is the number of shell boxes introduced in Step 3, and the matrix $B\colon\R^M\to\R^n$, constructed in Step 4, satisfies $BB^*\preceq I_n$. Let $K_t\coloneqq\{z:\norm{G_{Bz}}_\infty\le t\}$, a closed, convex and symmetric set (from \cref{sec:Sampling} on we abbreviate $G_z\coloneqq G_{Bz}$). Since $\eta=Bz$ has coordinate variances $(BB^*)_{ii}\le1$ and $n\le C_dM$, \eqref{eq:CoverIntro} gives $\gm(K_{2D})\ge\mathrm{e}^{-\kappa M}$ for a constant $\kappa$ depending only on $d$.

Whenever $\gm(K_t)\ge\mathrm{e}^{-\kappa M}$, the conditioned vector $Z$ with law $\nu\coloneqq\gm(\,\cdot\mid K_t)$ has differential entropy at least $\frac M2\log(2\pi)-\kappa M$, which is at most $\kappa+\frac12$ nats per coordinate below the entropy of $\gm$. Its Rosenblatt transform $U=(U_1,\dots,U_M)$, $U_j\coloneqq2F_j(Z_j\mid Z_{<j})-1$ \cite{Rosenblatt52}, built from the successive conditional distribution functions $F_j$, has independent coordinates that are uniform on $[-1,1]$, and an elementary inequality between entropy and the Gini mean difference, combined with the chain rule for entropy, yields
\begin{equation}\label{eq:EntropyIntro}
    \sum_{j=1}^M\E_{\nu}(Z_jU_j)\ge c_\kappa M
\end{equation}
with $c_\kappa>0$ depending only on $\kappa$ (\cref{prop:Correlation}). Each $\E_{\nu}(Z_jU_j)$ is nonnegative, and their average is bounded below by $c_\kappa$. The lower bound $\gm(K_t)\ge\mathrm{e}^{-\kappa M}$ enters the argument only through \eqref{eq:EntropyIntro}.

\paragraph{Step 3: test functions with linear response.}
The heart of the construction is a family of functions $W(u)$ on $\cube$, indexed by coefficient vectors $u\in[-1,1]^M$, with $\norm{W(u)}_1\le1$ for every $u$ and such that $\Jst\Phig W(u)$, the primitive of $W(u)$ after removal of its component in the \emph{low Haar cross}, the span of all tensor Haar functions of level below $h$, is linear in $u$ up to an error that is small in quadratic mean over $u$; see \eqref{eq:RemainderIntro}. The index set is a family $\Shell$ of $M$ dyadic boxes of volume $2^{-h}$ whose sides have depths in $3\N_0$ (hence $h\in3\N$), the \emph{shell boxes}, to each of which a localized Walsh function $q_R$ is attached; from now on we index the coordinates of $\R^M$ by $\Shell$. The test function is the normalized difference of two positive products:
\begin{equation}\label{eq:TestFunctionIntro}
    W(u)\coloneqq\tfrac12\bigl(Q_+(u)-Q_-(u)\bigr),\qquad Q_\pm(u)\coloneqq\prod_a\bigl(1\pm\varepsilon B_a(u)\bigr),\qquad\varepsilon=\tfrac14.
\end{equation}
Here the shapes of the shell boxes are grouped into blocks by the depth $a$ of their first side, and the block function $B_a(u)$ is obtained by inserting the layer functions $\sum_Ru_Rq_R$ (one for each shape of the block, $R$ running over the shell boxes of that shape) into the multilinear extension of $\xi\mapsto\tanh(\theta\sum_r\xi_r)$ from the Boolean cube (\cref{sec:Witness}). Because every $B_a$ is odd in one binary digit of the first coordinate that no block of smaller first depth sees, both products integrate to one along every line parallel to the first axis, which gives $\norm{W(u)}_1\le1$ (\cref{lem:Normalization}).

The linear part of $W(u)$ in $u$ is $\Wlin(u)=\sum_R\lambda_Ru_Rq_R$ with $\lambda_R\asymp\theta$, where the modulation parameter is $\theta\coloneqq t_d(h+1)^{-(d-2)/2}$ with a small constant $t_d$. The nonlinear part is not small as a function, but it becomes small after removal of its component in the low Haar cross and integration: \cref{thm:Remainder} shows that
\begin{equation}\label{eq:RemainderIntro}
    \E_{U}\paren[\big]{\norm{\Jst\Phig W(U)-\Jst\Wlin(U)}_2^2}
    \le C_d\theta^4\;\E_{U}\paren[\big]{\norm{\Jst\Wlin(U)}_2^2}
\end{equation}
for uniform random coefficients $U$, where $\Phig$ denotes the orthogonal projection onto the orthogonal complement of the low Haar cross. The relative root-mean-square error is thus $O_d(\theta^2)=O_d((h+1)^{-(d-2)})$, which tends to zero as $h\to\infty$ because $d\ge3$.

\paragraph{Step 4: sampling discretization.}
By \eqref{eq:DualityIntro}, $\inp{G_{Bz}}{W(u)}=z\cdot V(u)$ with the \emph{response} $V(u)\coloneqq B^*\bigl((\Jst W(u))(P_i)\bigr)_{i\le n}$. The points come from a sampling discretization theorem of Limonova and Temlyakov \cite[Theorem~1.1]{LimonovaTemlyakov22}, applied to a space of piecewise multiaffine functions that contains the primitives $\Jst q_R$ and the primitives of the low Haar cross, and the matrix $B$ is built from the quadrature weights, the evaluations of the primitives, a projection that annihilates the low Haar cross, and an inverse Gram matrix. Since the remainder $\Jst\Phig W(u)-\Jst\Wlin(u)$ need not lie in this space, the sampling space is enlarged by $\sqrt w$, where $w(p)\coloneqq\E_{U}(|(\Jst\Phig W(U)-\Jst\Wlin(U))(p)|^2)$ is the mean-square remainder; the quadrature rule then also bounds the weighted sum of the values $w(P_i)$ by a constant multiple of $\int w$, which transfers \eqref{eq:RemainderIntro} to the points (\cref{lem:SamplingLemma}).

\Cref{prop:Response} shows that the response is diagonal in $u$ up to a small error: there are numbers $s_R>0$, all comparable to a common scale $s\asymp\theta\,(h+1)^{(d-1)/2}2^{-h}$, which results from our chosen normalization of $B$, such that
\begin{equation}\label{eq:ResponseIntro}
    V(u)=\diag(s_R)\,u+e(u),\qquad s_*\coloneqq\min_Rs_R\asymp s,\qquad \E_{U}\paren[\big]{\norm{e(U)}_2^2}\le\rho_h^2s_*^2M,
\end{equation}
with $\rho_h\coloneqq C_\rho(h+1)^{-(d-2)}$ for a constant $C_\rho$ depending only on $d$.

\paragraph{Step 5: conclusion.}
For $z\in K_t$ and every $u\in[-1,1]^M$ we have $z\cdot V(u)=\inp{G_{Bz}}{W(u)}\le t$. Substituting $u=U(z)$, where $U$ is the Rosenblatt transform of $\nu$, taking expectations under $\nu$, and using \eqref{eq:ResponseIntro}, \eqref{eq:EntropyIntro} together with $\E_{\nu}(Z_RU_R)\ge0$ for every $R$, the exact uniform law of $U$ under $\nu$, and $\E_{\nu}(\norm{Z}_2^2)\le M$ (\cref{lem:SecondMoment}), we obtain
\begin{equation}\label{eq:ConclusionIntro}
    t\ \ge\ \E_{\nu}(Z\cdot V(U))\ \ge\ s_*\sum_R\E_{\nu}(Z_RU_R)-\sqrt M\,\bigl(\E_{U}\paren[\big]{\norm{e(U)}_2^2}\bigr)^{1/2}\ \ge\ (c_\kappa-\rho_h)\,s_*M
\end{equation}
(\cref{prop:Obstruction}). With $t=2D$ and $h$ so large that $\rho_h\le c_\kappa/2$, this gives $D\gtrsim s_*M$. The exponent $d-1/2$ arises as follows: the number $M\asymp(h+1)^{d-1}2^h$ of shell boxes and the factor $s/\theta\asymp(h+1)^{(d-1)/2}2^{-h}$ of the response scale give $s_*M\asymp\theta\,(h+1)^{3(d-1)/2}$, and the modulation parameter $\theta\asymp(h+1)^{-(d-2)/2}$ reduces this to $(h+1)^{d-1/2}\asymp(\log n)^{d-1/2}$.

\paragraph{Relation to earlier methods.}
Products of modulated block sums such as $Q_\pm$ in \eqref{eq:TestFunctionIntro} go back to Riesz products, which appear in Hal\'asz's lower bounds for the Lebesgue-measure discrepancy \cite{Halasz81} and in proofs of small-ball inequalities: in Temlyakov's Riesz-product proof \cite{Temlyakov95Haar} of Talagrand's two-dimensional Haar inequality \cite[Proposition~3.2]{Talagrand94} and its trigonometric analogue \cite{Temlyakov95}, and in the higher-dimensional results of Bilyk and Lacey \cite{BilykLacey08} and Bilyk, Lacey and Vagharshakyan \cite{BilykLaceyVagharshakyan08}, where coincidences between Haar functions of equal depth in the same coordinate are the main obstacle. Here the product is paired with the counting function of a Gaussian weight vector, its coefficients are the Rosenblatt coordinates of that vector, and the coincidences are controlled: since the depths of the shell boxes are spaced by three, every nonvanishing product of the functions $q_I$ in one coordinate reduces to a factor $\pm q_{I'}$ or $\pm\one_{I'}$, whose contributions are controlled by the smoothing of the block functions (\cref{lem:WalshTanh}) and by the estimate of the high-frequency primitives (\cref{sec:Remainder}). The passage from colorings to Gaussian vectors through hereditary discrepancy is in the spirit of the factorization-norm approach of \cite{MNT20}, but the quantity we estimate is a Gaussian small-ball probability rather than a $\gamma_2$ norm: hereditary discrepancy bounds it from below (\cref{lem:GaussianCover}), and the entropy argument bounds it from above for radii $t\ll(h+1)^{d-1/2}$ (\cref{prop:Obstruction}).

After the preliminaries in \cref{sec:Preliminaries}, Sections~\ref{sec:Witness}--\ref{sec:Entropy} develop the main construction and entropy argument, \cref{sec:Conclusion} completes the proof of \cref{thm:Main}, and \cref{sec:Consequences} establishes the extensions and further consequences.

\section{Preliminaries}
\label{sec:Preliminaries}

This section collects definitions, external results and auxiliary lemmas. \Cref{sec:HereditaryDiscrepancy} and the Gaussian inequalities of \cref{sec:GaussianTools} are used in Step~1 of the overview (\cref{sec:Conclusion}); the Rosenblatt transform and the entropy facts of \cref{sec:GaussianTools,sec:Entropy-Preliminaries} in Step~2 and the second-moment bound of \cref{lem:SecondMoment} in Step~5 (both in \cref{sec:Entropy}); the dyadic functions and the Walsh estimate of \cref{sec:Dyadic,sec:Walsh} in the construction of the test functions (\cref{sec:Witness,sec:Remainder}), the dyadic functions also in \cref{sec:Sampling}; and the sampling theorem of \cref{sec:SamplingPreliminaries} in \cref{sec:Sampling}. \Cref{lem:qProperties,lem:WalshTanh} are specific to our construction; the remaining lemmas are standard and included for completeness.

\subsection{Discrepancy and hereditary discrepancy}
\label{sec:HereditaryDiscrepancy}

For a real matrix $A$ with $n$ columns, its discrepancy and hereditary discrepancy are
\begin{equation*}
    \disc(A)\coloneqq\min_{\sigma\in\{-1,1\}^n}\norm{A\sigma}_\infty,
    \qquad
    \herdisc(A)\coloneqq\max_{S\subseteq[n]}\disc(A_S),
\end{equation*}
where $A_S$ is the submatrix of $A$ formed by the columns indexed by $S$. For a finite point set $P=\{P_1,\dots,P_n\}\subset\R^d$, the set system of anchored boxes with corners in the unit cube is $\{P\cap A(x):x\in\cube\}$; its incidence matrix $A_P\in\{0,1\}^{m\times n}$ has one row $(\one_{\{P_i\le x\}})_{i\le n}$ for each of the finitely many distinct sets $P\cap A(x)$. For every subset $S\subseteq P$, the columns of $A_P$ indexed by $S$ form, up to repeated rows, the incidence matrix of the anchored boxes on $S$, so that
\begin{equation}\label{eq:HerdiscPoints}
    \herdisc(A_P)=\max_{S\subseteq P}\ \min_{\chi\colon S\to\{-1,1\}}\ \sup_{x\in\cube}\ \Bigl|\sum_{p\in S\cap A(x)}\chi(p)\Bigr|
    \le\max_{S\subseteq P}\disc(S,\Ad).
\end{equation}
The inequality holds because $\Ad$ also contains anchored boxes with corners outside the unit cube; \eqref{eq:HerdiscPoints} is used in the proof of \cref{thm:Main}. We shall use the inequality between linear and hereditary discrepancy due to Lov\'asz, Spencer and Vesztergombi \cite[Corollary~1]{LSV86} (their discrepancy is defined with colorings in $\{0,1\}$ and is half of ours, so that their Corollary~1 is the first assertion of \cref{lem:Rounding}); see also \cite[Theorem~4.6]{Matousek99book}. We state it in the form \eqref{eq:RoundingCover}, which is all we need, and include the short proof by dyadic rounding for completeness.

\begin{lemma}[Rounding lemma]
    \label{lem:Rounding}
    Let $A\in\R^{m\times n}$ satisfy $\herdisc(A)\le D$. Then for every $y\in[0,1]^n$ there is $b\in\{0,1\}^n$ with $\norm{A(y-b)}_\infty\le D$. Consequently,
    \begin{equation}\label{eq:RoundingCover}
        A[-1,1]^n\subseteq\bigcup_{\sigma\in\{-1,1\}^n}\bigl(A\sigma+[-2D,2D]^{m}\bigr).
    \end{equation}
\end{lemma}

\begin{proof}
    Let $y$ first have dyadic coordinates with common denominator $2^k$, and iterate for $j=k,k-1,\dots,1$: let $S$ be the set of coordinates whose numerator (in denominator $2^j$) is odd, choose $\sigma\in\{-1,1\}^S$ with $\norm{A_S\sigma}_\infty\le D$, and add $2^{-j}\sigma$ to the coordinates in $S$. Odd numerators $o$ satisfy $0<o<2^j$, so all coordinates stay in $[0,1]$, all numerators become even, and $Ay$ changes by at most $2^{-j}D$ in the maximum norm. After the step $j=1$ the vector $b\in\{0,1\}^n$ satisfies $\norm{A(y-b)}_\infty\le D\sum_{j=1}^k2^{-j}\le D$. For general $y$, approximate by dyadic vectors $y^{(\ell)}\to y$, pass to a subsequence along which the resulting $b^{(\ell)}=b$ is constant, and let $\ell\to\infty$. For \eqref{eq:RoundingCover}, apply this to $y\coloneqq(\eta+\one)/2$ with $\eta\in[-1,1]^n$ and put $\sigma\coloneqq2b-\one$; then $A\eta-A\sigma=2A(y-b)$.
\end{proof}

\subsection{Dyadic intervals, Rademacher and Haar functions, primitives}
\label{sec:Dyadic}

For $j\in\N_0$, let $\dy_j$ be the set of the $2^j$ half-open dyadic intervals $[k2^{-j},(k+1)2^{-j})$, $0\le k<2^j$, of \emph{depth} $j$; we write $\depth(I)\coloneqq j$ for $I\in\dy_j$, and $|I|=2^{-j}$ for its length. Two dyadic intervals are either nested or disjoint. For $I\in\dy_j$ and $0\le b\le j$, the \emph{ancestor} of $I$ of depth $b$ is the unique interval of $\dy_b$ containing $I$. The \emph{Rademacher functions} are $\Rad_j(x)\coloneqq(-1)^{\lfloor2^jx\rfloor}$ for $x\in[0,1)$ and $j\ge1$; $\Rad_j$ is constant on every interval of $\dy_j$ and takes the values $+1$ and $-1$ on the left and right halves of every interval of $\dy_{j-1}$. The ($L^\infty$-normalized) \emph{Haar function} of $J\in\dy_b$ is
\begin{equation*}
    h_J\coloneqq\Rad_{b+1}\one_J=\one_{J^-}-\one_{J^+},
\end{equation*}
where $J^-$ and $J^+$ are the left and right halves of $J$. Its \emph{level} is $\lev(h_J)\coloneqq b$, and the constant function $\one$ has level $\lev(\one)\coloneqq-1$. The functions $\one$ and $h_J$, $J\in\bigcup_b\dy_b$, form an orthogonal basis of $L^2[0,1]$, and the \emph{tensor Haar functions}
\begin{equation*}
    \varphi_1\otimes\dots\otimes\varphi_d,\qquad\varphi_j\in\{\one\}\cup\{h_J:J\in\textstyle\bigcup_b\dy_b\},
\end{equation*}
form an orthogonal basis of $L^2(\cube)$. The level of a tensor Haar function is the sum $\sum_j\lev(\varphi_j)\in\{-d,-d+1,\dots\}$ of its coordinate levels; for every integer $L$ there are only finitely many tensor Haar functions of level at most $L$.

Recall from \eqref{eq:DualityIntro} the \emph{primitive}
\begin{equation}\label{eq:Primitive}
    (\Jst f)(p)\coloneqq\int_{[p,1]^d}f(x)\,dx,\qquad p\in\cube,
\end{equation}
of an integrable function $f$ on $\cube$. The function $\Jst f$ is continuous, and $\Jst$ is linear and acts on tensor products coordinatewise: $\Jst(f_1\otimes\dots\otimes f_d)=\Jst f_1\otimes\dots\otimes\Jst f_d$, where in one dimension $(\Jst g)(t)=\int_t^1g$ for $t\in[0,1]$. With the variable $t\in[0,1]$, one computes
\begin{equation}\label{eq:HaarPrimitives}
    \Jst\one=1-t,\quad\norm{\Jst\one}_2^2=\tfrac13,\qquad
    \Jst h_J=-\tau_J,\quad\norm{\Jst h_J}_2^2=\tfrac{|J|^3}{12},
\end{equation}
where $\tau_J$ is the tent function supported on $J$ that is affine on $J^-$ and on $J^+$ and has height $|J|/2$ at the midpoint of $J$.

\paragraph{Localized Walsh functions $q_I$ and their primitives $\psi_I$.}
The functions $q_I$ are the one-dimensional factors of the localized Walsh functions $q_R$ attached to the shell boxes in \cref{sec:Witness}. The three Rademacher factors make $q_I$ orthogonal to all polynomials of degree at most two, so that $\psi_I$ is orthogonal to affine functions (\cref{lem:qProperties}); this orthogonality is used in \cref{sec:Witness,sec:Sampling}. For $I\in\dy_a$ we define
\begin{equation}\label{eq:qDefinition}
    q_I\coloneqq\one_I\,\Rad_{a+1}\Rad_{a+2}\Rad_{a+3}
    =h_{I_0}-h_{I_1}-h_{I_2}+h_{I_3},
    \qquad
    \psi_I\coloneqq\Jst q_I,
\end{equation}
where $I_0,I_1,I_2,I_3\in\dy_{a+2}$ are the four quarters of $I$ from left to right. The second expression follows from the sign pattern of $\Rad_{a+1}\Rad_{a+2}$ on the quarters and from $\Rad_{a+3}\one_{I_k}=h_{I_k}$. Thus $|q_I|=\one_I$, $q_I$ is a combination of Haar functions of level $a+2$, and $q_I$ is constant on every interval of $\dy_{a+3}$. Up to sign, $q_I$ is one of the iterated Haar differences used in Kolmogorov's method for lower bounds on the error of cubature formulas in classes of functions with mixed smoothness; see \cite{Bilyk14}.

\begin{lemma}[Moments of $q_I$ and $\psi_I$]
    \label{lem:qProperties}
    Let $I\in\dy_a$. Then $\int q_I(x)x^k\,dx=0$ for $k=0,1,2$. The primitive $\psi_I$ is continuous, vanishes outside $I$, is affine on every interval of $\dy_{a+3}$, and satisfies
    \begin{equation}\label{eq:psiProperties}
        \int_0^1\psi_I(p)\,dp=\int_0^1p\,\psi_I(p)\,dp=0,\qquad\norm{\psi_I}_2^2=\frac{|I|^3}{192}.
    \end{equation}
    In particular, $\inp{\psi_I}{g}=0$ for every function $g$ that is affine on an interval containing $I$.
\end{lemma}

\begin{proof}
    Write $I=[c,c+\ell)$ and substitute $x=c+\ell u$ with $u\in[0,1)$; then $\Rad_{a+k}(x)=\Rad_k(u)$ for $k=1,2,3$, so the moments are $\ell\int_0^1\Rad_1\Rad_2\Rad_3(u)\,(c+\ell u)^k\,du$. Writing $u=\sum_{i\ge1}\iota_i2^{-i}$ with binary digits $\iota_i=(1-\Rad_i(u))/2$, the function $u^k$ with $k\le2$ is an absolutely convergent series of products of at most two Rademacher functions, while $\Rad_1\Rad_2\Rad_3$ is a product of three. Since the products $\prod_{i\in T}\Rad_i$ over distinct finite sets $T\subset\N$ are orthonormal in $L^2[0,1]$, integrating termwise shows that the moments vanish for $k\le2$. Since $q_I$ vanishes outside $I$ and $\int q_I=0$, the primitive $\psi_I(p)=\int_p^1q_I$ vanishes for $p\notin I$. Fubini's theorem gives $\int_0^1\psi_I(p)\,dp=\int q_I(x)x\,dx=0$ and $\int_0^1p\,\psi_I(p)\,dp=\frac12\int q_I(x)x^2\,dx=0$; as $\psi_I$ vanishes outside $I$, the assertion on affine functions $g$ follows. On $I$, the function $q_I$ equals $\pm1$ on the eight intervals of $\dy_{a+3}$ contained in $I$, with the sign pattern $(+,-,-,+,-,+,+,-)$; integrating from the right, $\psi_I$ is the sum of four tents of height $\ell/8$ and base width $\ell/4$ with signs $(-,+,+,-)$, which gives $\norm{\psi_I}_2^2=4\cdot\frac{\ell}{4}\cdot\frac{(\ell/8)^2}{3}=\frac{\ell^3}{192}$. The remaining assertions (continuity, affinity on the intervals of $\dy_{a+3}$) follow from this representation as a sum of tents.
\end{proof}

\subsection{Fourier--Walsh expansions}
\label{sec:Walsh}

Every function $F\colon\{-1,1\}^m\to\R$ has the expansion $F=\sum_{T\subseteq[m]}\widehat F(T)\chi_T$ with $\chi_T(\xi)\coloneqq\prod_{i\in T}\xi_i$ and $\widehat F(T)=\E(F(\xi)\chi_T(\xi))$ for a uniformly distributed $\xi$; see \cite[Chapter~1 and Section~2.2]{ODonnell14}. The characters $\chi_T$ are orthonormal, so $\E(F(\xi)^2)=\sum_T\widehat F(T)^2$. The \emph{multilinear extension} of $F$ is
\begin{equation}\label{eq:MultilinearExtension}
    \widetilde F(y)\coloneqq\sum_{T\subseteq[m]}\widehat F(T)\prod_{i\in T}y_i,\qquad y\in[-1,1]^m.
\end{equation}
If $\xi_1,\dots,\xi_m$ are independent random signs with $\E(\xi_i)=y_i$, then $\widetilde F(y)=\E(F(\xi))$; in particular $|\widetilde F(y)|\le\max|F|$. For $i\in[m]$, the discrete derivative $D_iF(\xi)\coloneqq\frac12\bigl(F(\xi^{i\to+})-F(\xi^{i\to-})\bigr)$, where $\xi^{i\to\pm}$ agrees with $\xi$ except that its $i$-th coordinate is $\pm1$, has the expansion $D_iF=\sum_{T\ni i}\widehat F(T)\chi_{T\setminus\{i\}}$ \cite[Proposition~2.19]{ODonnell14}. The operators $D_i$ commute, and for distinct $i,j,k$,
\begin{equation}\label{eq:DerivativeParseval}
    \E\paren[\big]{(D_iF)^2}=\sum_{T\ni i}\widehat F(T)^2,
    \qquad
    \E\paren[\big]{(D_iD_jD_kF)^2}=\sum_{T\supseteq\{i,j,k\}}\widehat F(T)^2.
\end{equation}

The block functions of \cref{sec:Witness} are multilinear extensions of the function in the following lemma. Part~(ii) identifies their linear part, and part~(iii) controls the nonlinear part in \cref{sec:Remainder}.

\begin{lemma}[Walsh coefficients of a smoothed block sum]
    \label{lem:WalshTanh}
    Let $m\in\N$, $\theta>0$, and $F(\xi)\coloneqq\tanh\bigl(\theta\sum_{i=1}^m\xi_i\bigr)$ for $\xi\in\{-1,1\}^m$. Then:
    \begin{enumerate}[label=\textup{(\roman*)}]
        \item $\widehat F(T)=0$ whenever $|T|$ is even;
        \item the first-order coefficient $b\coloneqq\widehat F(\{i\})$ does not depend on $i$, and $\theta\bigl(1-\theta^2(m-\tfrac23)\bigr)\le b\le\theta$;
        \item for every $\calE\subseteq[m]$,
        \begin{equation}\label{eq:WalshActivity}
            \sum_{\emptyset\ne T\subseteq\calE}\widehat F(T)^2\le\theta^2|\calE|
            \qquad\text{and}\qquad
            \sum_{T\subseteq\calE,\ |T|\ge3}\widehat F(T)^2\le\theta^6|\calE|^3.
        \end{equation}
    \end{enumerate}
\end{lemma}

\begin{proof}
    (i) $F$ is odd, $F(-\xi)=-F(\xi)$, and $\chi_T(-\xi)=(-1)^{|T|}\chi_T(\xi)$, so $\widehat F(T)=-(-1)^{|T|}\widehat F(T)$.

    (ii) By symmetry $b$ does not depend on $i$. Let $S\coloneqq\sum_{l\ne i}\xi_l$ and $g(s)\coloneqq\tanh(\theta s)$. Then $b=\E(\xi_iF(\xi))=\frac12\E(g(S+1)-g(S-1))=\frac12\E(\int_{-\theta}^{\theta}\sech^2(\theta S+t)\,dt)$. Since $\sech^2\le1$, $b\le\theta$. Since $\sech^2(y)\ge1-y^2$ and $\E((\theta S+t)^2)=\theta^2(m-1)+t^2$,
    \begin{equation*}
        b\ge\frac12\int_{-\theta}^{\theta}\bigl(1-\theta^2(m-1)-t^2\bigr)dt=\theta-\theta^3(m-1)-\frac{\theta^3}{3}.
    \end{equation*}

    (iii) Since $\tanh$ is $1$-Lipschitz, $|D_iF|\le\theta$, and \eqref{eq:DerivativeParseval} gives $\sum_{\emptyset\ne T\subseteq\calE}\widehat F(T)^2\le\sum_{i\in\calE}\sum_{T\ni i}\widehat F(T)^2\le\theta^2|\calE|$. For distinct $i,j,k$, with $g$ as in (ii) and $S'\coloneqq\sum_{l\notin\{i,j,k\}}\xi_l$,
    \begin{equation*}
        D_iD_jD_kF(\xi)=\frac18\bigl(g(S'+3)-3g(S'+1)+3g(S'-1)-g(S'-3)\bigr)=g'''(\zeta)
    \end{equation*}
    for some $\zeta\in(S'-3,S'+3)$, by the mean value theorem for third-order differences with step two. As $|\tanh'''|\le2$, we get $|D_iD_jD_kF|\le2\theta^3$. Every $T\subseteq\calE$ with $|T|\ge3$ contains a triple $\{i,j,k\}\subseteq\calE$, so by \eqref{eq:DerivativeParseval}
    \begin{equation*}
        \sum_{T\subseteq\calE,\ |T|\ge3}\widehat F(T)^2\le\sum_{\{i,j,k\}\subseteq\calE}\E\paren[\big]{(D_iD_jD_kF)^2}\le\binom{|\calE|}{3}4\theta^6\le\theta^6|\calE|^3.\qedhere
    \end{equation*}
\end{proof}

\subsection{Gaussian inequalities and the Rosenblatt transform}
\label{sec:GaussianTools}

We call a random vector $\eta$ in $\R^n$ a \emph{centered Gaussian vector} if $\eta=Bz$ for some $k\in\N$, a standard Gaussian vector $z$ in $\R^k$ and a matrix $B\in\R^{n\times k}$; its covariance $BB^*$ ($B^*$ the transpose) may be singular. We use two classical inequalities, a second-moment bound and the Rosenblatt transform.

\begin{literaturetheorem}[{\v{S}id\'ak \cite[Corollary~1]{Sidak67}}]
    \label{thm:Sidak}
    Let $\eta$ be a centered Gaussian vector in $\R^n$ and $c_1,\dots,c_n>0$. Then
    \begin{equation*}
        \Prob(|\eta_i|\le c_i\ \text{for all }i\le n)\ge\prod_{i=1}^n\Prob(|\eta_i|\le c_i).
    \end{equation*}
\end{literaturetheorem}

\begin{literaturetheorem}[{Anderson \cite[Theorem~1 and Corollary~2]{Anderson55}}]
    \label{thm:Anderson}
    Let $\eta$ be a centered Gaussian vector in $\R^m$ and let $K\subset\R^m$ be closed, convex and symmetric ($K=-K$). Then $\Prob(\eta\in K+c)\le\Prob(\eta\in K)$ for every $c\in\R^m$.
\end{literaturetheorem}

Both inequalities are applied to possibly singular covariances. \v{S}id\'ak proves his inequality for zero means, arbitrary variances and an arbitrary correlation matrix and obtains the singular case by approximation with nonsingular distributions \cite[Theorem~1, Corollary~1 and p.~628]{Sidak67}; coordinates with zero variance are trivial. Anderson proved his inequality for symmetric unimodal densities \cite[Theorem~1 and Corollary~2]{Anderson55}, which covers nondegenerate centered Gaussian vectors (see also \cite[Section~1.8]{Bogachev98}); the degenerate case follows by adding an independent noise $\delta g$ with $g$ standard Gaussian, applying the nondegenerate case to the enlarged set $K+[-\rho,\rho]^m$, and letting first $\delta\to0$ and then $\rho\to0$, using that $K$ is closed.

\begin{lemma}[Second moment of a conditioned Gaussian vector]
    \label{lem:SecondMoment}
    Let $K\subset\R^M$ be closed, convex and symmetric with $\gm(K)>0$, and let $Z$ be standard Gaussian in $\R^M$. Then $\E(\norm{Z}_2^2\mid Z\in K)\le M$.
\end{lemma}

\begin{proof}
    Write $Z=R\Theta$ with $R\coloneqq\norm{Z}_2$ and $\Theta\coloneqq Z/\norm{Z}_2$; then $R$ and $\Theta$ are independent and $\Theta$ is uniformly distributed on the unit sphere. Let $a(r)\coloneqq\Prob(r\Theta\in K)$. Since $K$ is convex, symmetric and nonempty, $0\in K$, so $r\omega\in K$ implies $r'\omega\in K$ for $0\le r'\le r$ and every unit vector $\omega$; hence $a$ is nonincreasing. By independence, $\E(\norm{Z}_2^2\one_{\{Z\in K\}})=\E(R^2a(R))$, and Chebyshev's association inequality, applied to the nondecreasing function $r\mapsto r^2$ and the nonincreasing function $a$ of the random variable $R$, gives $\E(R^2a(R))\le\E(R^2)\E(a(R))=M\gm(K)$.
\end{proof}

For a random vector $Z=(Z_1,\dots,Z_M)$ and $j\le M$ we write $Z_{<j}\coloneqq(Z_1,\dots,Z_{j-1})$ and $Z_{\le j}\coloneqq(Z_1,\dots,Z_j)$, and similarly $z_{<j}$ and $z_{\le j}$ for points $z\in\R^M$.

\begin{lemma}[Rosenblatt transform \cite{Rosenblatt52}]
    \label{lem:Rosenblatt}
    Let $Z=(Z_1,\dots,Z_M)$ be a random vector with a density on $\R^M$. For $j\le M$ let $F_j(\,\cdot\mid z_{<j})$ be a regular version of the conditional distribution function of $Z_j$ given $Z_{<j}=z_{<j}$ (the distribution function of $Z_1$ for $j=1$), and put $U_j\coloneqq2F_j(Z_j\mid Z_{<j})-1$. Then $U_1,\dots,U_M$ are independent and uniformly distributed on $[-1,1]$.
\end{lemma}

\begin{proof}
    For $\Prob_{Z_{<j}}$-almost every $z_{<j}$, where $\Prob_{Z_{<j}}$ is the law of $Z_{<j}$, the conditional law of $Z_j$ given $Z_{<j}=z_{<j}$ has a density, so that $F_j(\,\cdot\mid z_{<j})$ is continuous and $U_j$ is uniform on $[-1,1]$ conditionally on $Z_{<j}$; hence $U_j$ is uniform on $[-1,1]$ and independent of $U_1,\dots,U_{j-1}$, which are functions of $Z_{<j}$. Induction on $j$ gives the claim.
\end{proof}

\subsection{Differential entropy}
\label{sec:Entropy-Preliminaries}

For a random vector $Z$ with density $\varphi$ on $\R^k$, the differential entropy is
\begin{equation*}
    \ent(Z)\coloneqq-\int\varphi\log\varphi,
\end{equation*}
whenever this integral exists; we also write $\ent(\mu)\coloneqq\ent(Z)$ if $Z$ has law $\mu$. If $\varphi$ is bounded with bounded support, then $\varphi\log\varphi$ is bounded and integrable, so $\ent(Z)$ is finite; moreover all marginal densities of $Z$ are bounded with bounded support. We use two standard facts.

For a random vector $Z=(Z_1,\dots,Z_M)$ with a bounded density with bounded support, let $\varphi_{\le j}$ and $\varphi_{<j}$ be the densities of $Z_{\le j}$ and $Z_{<j}$ (with $\varphi_{<1}\coloneqq1$), and let $\Prob_{Z_{<j}}$ be the law of $Z_{<j}$. For $\Prob_{Z_{<j}}$-almost every $z_{<j}$, namely whenever $\varphi_{<j}(z_{<j})>0$, the conditional law of $Z_j$ given $Z_{<j}=z_{<j}$ has the bounded density $\varphi_{\le j}(z_{<j},\cdot)/\varphi_{<j}(z_{<j})$ with bounded support; we write $\ent(Z_j\mid Z_{<j}=z_{<j})$ for its entropy, which is finite, and $\ent(Z_j\mid Z_{<j})\coloneqq\E(\ent(Z_j\mid Z_{<j}=Z_{<j}))$ for its average, so that $\ent(Z_1\mid Z_{<1})=\ent(Z_1)$.

\begin{lemma}[Chain rule]
    \label{lem:ChainRule}
    Let $Z=(Z_1,\dots,Z_M)$ have a bounded density with bounded support. Then each $\ent(Z_j\mid Z_{<j})$ is finite and
    \begin{equation*}
        \ent(Z)=\sum_{j=1}^M\ent(Z_j\mid Z_{<j}).
    \end{equation*}
\end{lemma}

\begin{proof}
    For $z_{<j}$ with $\varphi_{<j}(z_{<j})>0$, the conditional density is $\varphi_{\le j}(z_{<j},\cdot)/\varphi_{<j}(z_{<j})$, and
    \begin{equation*}
        \ent(Z_j\mid Z_{<j})=-\int\varphi_{\le j}\log\varphi_{\le j}+\int\varphi_{\le j}\log\varphi_{<j}=\ent(Z_{\le j})-\ent(Z_{<j}),
    \end{equation*}
    where Fubini's theorem applies because $\varphi_{\le j}|\log\varphi_{\le j}|$ and $\varphi_{<j}|\log\varphi_{<j}|$ are integrable. With $\ent(Z_{<1})\coloneqq0$, summing over $j$ gives the claim; see also \cite[Theorem~8.6.2]{CoverThomas06}.
\end{proof}

\begin{lemma}[Entropy and mean absolute deviation]
    \label{lem:MaxEntropy}
    Let $Z$ be a real random variable with a density, finite entropy and $\E(|Z|)<\infty$, and let $\beta\coloneqq\E(|Z-\E Z|)$, which is positive since $Z$ has a density. Then $\ent(Z)\le\log(2\mathrm{e}\beta)$.
\end{lemma}

\begin{proof}
    Let $\varphi$ be the density of $Z$, $\bar z\coloneqq\E Z$, and $\ell(t)\coloneqq(2\beta)^{-1}\mathrm{e}^{-|t-\bar z|/\beta}$. Gibbs' inequality $\int\varphi\log(\varphi/\ell)\ge0$ gives $\ent(Z)\le-\int\varphi\log\ell=\log(2\beta)+\beta^{-1}\E(|Z-\bar z|)=\log(2\mathrm{e}\beta)$.
\end{proof}

\subsection{Sampling discretization}
\label{sec:SamplingPreliminaries}

The points of our construction are provided by a sampling discretization theorem of Limonova and Temlyakov. We state the special case that we need: it is their Theorem~1.1 under their Condition~E with parameter $t=1$, which is the pointwise bound $\sum_k|u_k(x)|^2\le D$ below.

\begin{literaturetheorem}[{Limonova and Temlyakov \cite[Theorem~1.1]{LimonovaTemlyakov22}}]
    \label{thm:Sampling}
    There are absolute constants $C_1,C_2,C_3>0$ with the following property. Let $\mu$ be a probability measure on a nonempty set $\Omega\subset\R^d$, and let $u_1,\dots,u_D$ be real functions on $\Omega$ that are orthonormal in $L^2(\Omega,\mu)$ and satisfy $\sum_{k=1}^D|u_k(x)|^2\le D$ for every $x\in\Omega$. Then there are points $\xi^1,\dots,\xi^m\in\Omega$ with $m\le C_1D$ such that
    \begin{equation*}
        C_2\norm{f}_{L^2(\mu)}^2\le\frac1m\sum_{j=1}^m|f(\xi^j)|^2\le C_3\norm{f}_{L^2(\mu)}^2
        \qquad\text{for all }f\in\operatorname{span}\{u_1,\dots,u_D\}.
    \end{equation*}
\end{literaturetheorem}

The points $\xi^j$ need not be distinct; distinct points are obtained in \cref{lem:SamplingLemma}. The lower inequality implies $m\ge D$, since the evaluation map $f\mapsto(f(\xi^j))_{j\le m}$ is injective on the $D$-dimensional space $\operatorname{span}\{u_1,\dots,u_D\}$; this is used in the proof of \cref{lem:SamplingLemma}.

\section{The test functions}
\label{sec:Witness}

We construct the test functions $W(u)$, $u\in[-1,1]^M$, of Step~3 in the overview, prove that they have $L^1$-norm at most one, and identify their linear part. The construction is parametrized by a depth $h\in3\N$, and we write
\begin{equation}\label{eq:vDefinition}
    v\coloneqq2^{-h}.
\end{equation}

\subsection{Shell boxes and localized Walsh functions}
\label{sec:ShellBoxes}

A \emph{shape} is a vector $r\in(3\N_0)^d$ with $|r|_1=h$; the set of shapes is $\Shapes$, and its cardinality is
\begin{equation}\label{eq:ShapeCount}
    K_h\coloneqq|\Shapes|=\binom{h/3+d-1}{d-1},\qquad c_d(h+1)^{d-1}\le K_h\le C_d(h+1)^{d-1}.
\end{equation}
For $r\in\Shapes$, the \emph{shell boxes of shape $r$} are the dyadic boxes $R=R_1\times\dots\times R_d$, that is, products of dyadic intervals, with $R_j\in\dy_{r_j}$; there are $2^h=1/v$ of them, each of volume $v$. We write $\Shell$ for the set of all shell boxes, $\shape(R)$ for the shape of $R\in\Shell$, and
\begin{equation}\label{eq:MDefinition}
    M\coloneqq|\Shell|=K_h/v.
\end{equation}
To every $R\in\Shell$ we attach the \emph{localized Walsh function} $q_R$ and its normalized primitive $T_R$,
\begin{equation}\label{eq:qRTR}
    q_R\coloneqq q_{R_1}\otimes\dots\otimes q_{R_d},
    \qquad
    T_R\coloneqq v^{-1}\,\psi_{R_1}\otimes\dots\otimes\psi_{R_d},
\end{equation}
with $q_I$ and $\psi_I$ as in \eqref{eq:qDefinition}. Then $|q_R|=\one_R$, so $\norm{q_R}_1=v$ and $\norm{q_R}_2^2=v$. Moreover, $q_R$ is a combination of tensor Haar functions of level $\sum_j(r_j+2)=h+2d$, and $\Jst q_R=vT_R$.

\begin{lemma}[Orthogonality of the primitives]
    \label{lem:TOrthogonality}
    For all $R,S\in\Shell$, $\inp{T_R}{T_S}=\gamma_dv\,\one_{\{R=S\}}$, where $\gamma_d\coloneqq192^{-d}$.
\end{lemma}

\begin{proof}
    We have $\inp{T_R}{T_S}=v^{-2}\prod_{j=1}^d\inp{\psi_{R_j}}{\psi_{S_j}}$. If $R=S$, then \eqref{eq:psiProperties} gives $v^{-2}\prod_j|R_j|^3/192=v^{-2}\,v^3\,192^{-d}=\gamma_dv$. If $R\ne S$ have the same shape, then $R_j\ne S_j$ for some $j$, these intervals have the same depth and are therefore disjoint, and $\inp{\psi_{R_j}}{\psi_{S_j}}=0$. If the shapes differ, there is a coordinate $j$ with $r_j\ne s_j$, say $r_j<s_j$; then $s_j\ge r_j+3$ because the shapes have coordinates in $3\N_0$. Either $S_j\cap R_j=\emptyset$, or $S_j$ is contained in an interval of $\dy_{r_j+3}$ inside $R_j$, on which $\psi_{R_j}$ is affine by \cref{lem:qProperties}; in both cases $\inp{\psi_{R_j}}{\psi_{S_j}}=0$ by the last assertion of \cref{lem:qProperties}.
\end{proof}

\subsection{Blocks and the positive product}
\label{sec:PositiveProduct}

The shapes are grouped according to their first coordinate. For $a\in3\N_0$ with $a\le h$, the \emph{block} $\Bl a\coloneqq\{r\in\Shapes:r_1=a\}$ has cardinality
\begin{equation}\label{eq:BlockSize}
    m_a\coloneqq|\Bl a|=\binom{(h-a)/3+d-2}{d-2}\le(h+1)^{d-2},
\end{equation}
where the inequality follows from $\binom{n+k}{k}=\prod_{i=1}^k\frac{n+i}{i}\le(n+1)^k$ and $(h-a)/3+1\le h+1$. There are $L_h\coloneqq h/3+1$ blocks, indexed by $a\in\{0,3,\dots,h\}$; we also refer to the index $a$ itself as a block. We fix the parameters
\begin{equation}\label{eq:Parameters}
    \varepsilon\coloneqq\tfrac14,\qquad
    \theta\coloneqq t_d\,(h+1)^{-(d-2)/2},\qquad
    \lambda\coloneqq\varepsilon\theta,
\end{equation}
where $t_d\in(0,\tfrac12]$ is a constant depending only on $d$; it is chosen in \cref{sec:Remainder}, see \eqref{eq:tdChoice}. By \eqref{eq:BlockSize},
\begin{equation}\label{eq:tdFirstCondition}
    \theta^2m_a\le t_d^2\le\tfrac14\qquad\text{for every block } a.
\end{equation}
Thus the sign sums $\theta\sum_{r\in\Bl a}\xi_r$ defining the functions $\Phi_a$ below have variance at most $\frac14$; \cref{lem:WalshTanh} controls the resulting nonlinear Fourier--Walsh terms.

A \emph{coefficient vector} is a vector $u=(u_R)_{R\in\Shell}\in[-1,1]^{\Shell}$; we fix an enumeration of $\Shell$ and identify $[-1,1]^{\Shell}$ with $[-1,1]^M$ and $\R^{\Shell}$ with $\R^M$. For a shape $r$, the \emph{layer function}
\begin{equation}\label{eq:LayerFunction}
    f_r(u;x)\coloneqq\sum_{R\in\Shell:\ \shape(R)=r}u_R\,q_R(x),\qquad x\in\cube,
\end{equation}
satisfies $|f_r(u;x)|\le1$, because the shell boxes of shape $r$ are disjoint and $|q_R|=\one_R$. For a block $a$, let $\Phi_a\colon\{-1,1\}^{\Bl a}\to[-1,1]$, $\Phi_a(\xi)\coloneqq\tanh\bigl(\theta\sum_{r\in\Bl a}\xi_r\bigr)$, and let $\widetilde\Phi_a$ be its multilinear extension \eqref{eq:MultilinearExtension}. The \emph{block function} is
\begin{equation}\label{eq:BlockFunction}
    B_a(u;x)\coloneqq\widetilde\Phi_a\bigl((f_r(u;x))_{r\in\Bl a}\bigr)
    =\sum_{T\subseteq\Bl a}\widehat\Phi_a(T)\prod_{r\in T}f_r(u;x),
\end{equation}
and $|B_a(u;x)|\le1$ by the remark following \eqref{eq:MultilinearExtension}. Finally, we define the \emph{positive products} and the \emph{test function}
\begin{equation}\label{eq:TestFunction}
    Q_\pm(u;x)\coloneqq\prod_{a}\bigl(1\pm\varepsilon B_a(u;x)\bigr),
    \qquad
    W(u;x)\coloneqq\tfrac12\bigl(Q_+(u;x)-Q_-(u;x)\bigr),
\end{equation}
where the products run over all $L_h$ blocks. Since $\varepsilon|B_a|\le1/4$, we have $Q_\pm\ge(3/4)^{L_h}>0$. Expanding the products,
\begin{equation}\label{eq:TestFunctionExpansion}
    W(u;x)=\sum_{\substack{S\subseteq\{0,3,\dots,h\}\\ |S|\text{ odd}}}\varepsilon^{|S|}\prod_{a\in S}B_a(u;x).
\end{equation}

\subsection{Normalization}
\label{sec:Normalization}

The following lemma is the reason for the product structure of $W$ and for grouping the shapes by their first coordinate: integrating along the first coordinate removes one block at a time.

\begin{lemma}[Normalization]
    \label{lem:Normalization}
    For every $u\in[-1,1]^{\Shell}$ and every $x'=(x_2,\dots,x_d)\in[0,1)^{d-1}$,
    \begin{equation*}
        \int_0^1Q_\pm(u;x_1,x')\,dx_1=1.
    \end{equation*}
    Consequently, $\norm{Q_\pm(u)}_1=1$ and $\norm{W(u)}_1\le1$ for every $u\in[-1,1]^{\Shell}$.
\end{lemma}

\begin{proof}
    Fix $u$ and $x'$, and consider all functions as functions of $x_1\in[0,1)$. Let $a$ be a block and $r\in\Bl a$. For $x_1\in I\in\dy_a$, the unique shell box $R$ of shape $r$ containing $(x_1,x')$ is $R=I\times R_2\times\dots\times R_d$, where $R_j\in\dy_{r_j}$ is the interval containing $x_j$; these $R_j$ do not depend on $x_1$. Hence
    \begin{equation*}
        f_r(u;x_1,x')=q_I(x_1)\,c_r(I),\qquad c_r(I)\coloneqq u_R\prod_{j=2}^dq_{R_j}(x_j)\in[-1,1],\qquad x_1\in I,
    \end{equation*}
    and $c_r(I)$ depends on $x_1$ only through $I$. By \cref{lem:WalshTanh}(i), only sets $T$ of odd cardinality occur in \eqref{eq:BlockFunction}, and $q_I^{|T|}=q_I$ on $I$ for odd $|T|$. Therefore, for $x_1\in I$,
    \begin{equation*}
        B_a(u;x_1,x')=q_I(x_1)\,\zeta_a(I),\qquad\zeta_a(I)\coloneqq\sum_{T\subseteq\Bl a}\widehat\Phi_a(T)\prod_{r\in T}c_r(I)=\widetilde\Phi_a\bigl((c_r(I))_{r\in\Bl a}\bigr)\in[-1,1].
    \end{equation*}
    Since $q_I=\one_I\Rad_{a+1}\Rad_{a+2}\Rad_{a+3}$, we conclude that
    \begin{equation}\label{eq:BlockStructure}
        B_a(u;\cdot,x')=\Rad_{a+3}\,g_a,\qquad g_a\coloneqq\Rad_{a+1}\Rad_{a+2}\sum_{I\in\dy_a}\zeta_a(I)\one_I,
    \end{equation}
    where $g_a$ is constant on every interval of $\dy_{a+2}$. In particular, $B_a(u;\cdot,x')$ is constant on every interval of $\dy_{a+3}$, and $\int_JB_a(u;x_1,x')\,dx_1=0$ for every $J\in\dy_{a+2}$, because $\Rad_{a+3}$ has integral zero over $J$.

    We remove the factors of $Q_\pm$ one at a time. Let $a$ be the largest block index among the remaining factors. Every other remaining block $a'$ satisfies $a'+3\le a$, so its factor $1\pm\varepsilon B_{a'}$ is constant on the intervals of $\dy_{a}$ and hence on those of $\dy_{a+2}$, while $\int_J(1\pm\varepsilon B_a)\,dx_1=|J|$ for every $J\in\dy_{a+2}$. Integrating over each $J\in\dy_{a+2}$ and summing, we see that the integral over $[0,1)$ of the product of the remaining factors does not change when the factor $1\pm\varepsilon B_a$ is removed. Continuing with the remaining factors until the empty product $1$ is reached gives $\int_0^1Q_\pm(u;x_1,x')\,dx_1=1$. As $Q_\pm>0$, Fubini's theorem yields $\norm{Q_\pm(u)}_1=1$, and $\norm{W(u)}_1\le\frac12(\norm{Q_+(u)}_1+\norm{Q_-(u)}_1)=1$.
\end{proof}

\subsection{The linear part}
\label{sec:LinearPart}

Each $B_a(u;x)$ is a polynomial in the coefficients $u$, and so is $W(u;x)$; we isolate the part of $W(u;x)$ of degree one in $u$. By \cref{lem:WalshTanh}(ii) with $m=m_a$, the first-order Fourier--Walsh coefficient $b_a\coloneqq\widehat\Phi_a(\{r\})$, $r\in\Bl a$, does not depend on $r$. Put
\begin{equation}\label{eq:LinearCoefficients}
    \lambda_R\coloneqq\varepsilon\,b_{a}\qquad\text{for } R\in\Shell \text{ with } \shape(R)\in\Bl a,
\end{equation}
and define
\begin{equation}\label{eq:WlinFlin}
    \Wlin(u)\coloneqq\sum_{R\in\Shell}\lambda_Ru_Rq_R,
    \qquad
    \Flin(u)\coloneqq\Jst\Wlin(u)=v\sum_{R\in\Shell}\lambda_Ru_RT_R,
\end{equation}
where $\Flin=\Jst\Wlin$ follows from $\Jst q_R=vT_R$.

\begin{lemma}[Linear part]
    \label{lem:LinearPart}
    For every block $a$, $\theta/2\le b_a\le\theta$; hence $\frac12\lambda\le\lambda_R\le\lambda$ for every $R\in\Shell$. The difference $W(u)-\Wlin(u)$ is a polynomial in $u$ all of whose monomials have degree at least three. Moreover, if $U$ is uniformly distributed on $[-1,1]^{\Shell}$, then
    \begin{equation}\label{eq:FlinNorm}
    \begin{gathered}
        \E_{U}\paren[\big]{\norm{\Flin(U)}_2^2}=\frac{\gamma_dv^3}{3}\sum_{R\in\Shell}\lambda_R^2,\\
        c_d\lambda^2(h+1)^{d-1}v^2\le\E_{U}\paren[\big]{\norm{\Flin(U)}_2^2}\le C_d\lambda^2(h+1)^{d-1}v^2.
    \end{gathered}
    \end{equation}
\end{lemma}

\begin{proof}
    \Cref{lem:WalshTanh}(ii) with $m=m_a$ gives $b_a\le\theta$ and $b_a\ge\theta\bigl(1-\theta^2(m_a-\tfrac23)\bigr)\ge\theta(1-\theta^2m_a)\ge\theta/2$ by \eqref{eq:tdFirstCondition}. The layer functions $f_r(u)$ are linear in $u$, and distinct layer functions involve disjoint sets of coefficients. By \eqref{eq:BlockFunction} and \cref{lem:WalshTanh}(i), the part of $B_a(u)$ of degree one in $u$ is $b_a\sum_{r\in\Bl a}f_r(u)$, and all other monomials of $B_a(u)$ have odd degree at least three. In \eqref{eq:TestFunctionExpansion}, a monomial of degree one can only arise from a singleton $S=\{a\}$ and the degree-one part of $B_a$; this gives $\sum_a\varepsilon b_a\sum_{r\in\Bl a}f_r(u)=\Wlin(u)$. All remaining monomials have degree at least three. Finally, since $\E(U_RU_S)=\frac13\one_{\{R=S\}}$, \cref{lem:TOrthogonality} gives $\E_{U}(\norm{\Flin(U)}_2^2)=v^2\sum_R\lambda_R^2\cdot\frac13\gamma_dv$, and \eqref{eq:FlinNorm} follows from $\lambda/2\le\lambda_R\le\lambda$, $M=K_h/v$ and \eqref{eq:ShapeCount}.
\end{proof}

\begin{remark}[Depth separation and smoothing]
    \label{rem:Construction}
    The depths of the shell boxes are multiples of three, so that two different depths in the same coordinate differ by at least three. Since $q_I$ is constant on every dyadic interval of depth at least $\depth(I)+3$, every product of functions $q_I$ in one coordinate is therefore either zero or, up to sign, equal to $q_{I'}$ or $\one_{I'}$ for the deepest interval $I'$ involved, according to the parity of the number of factors of that depth (\cref{lem:CellMonomial}); likewise, the proof of \cref{lem:Normalization} can peel off one block at a time. The smoothing $\tanh$ in the block functions keeps $|B_a|\le1$, and hence $1\pm\varepsilon B_a\ge\frac34$, for all coefficient vectors and all points, while by \cref{lem:WalshTanh} the function $\Phi_a$ is linear up to cubic and higher Fourier--Walsh coefficients whose squares sum to at most $\theta^6m_a^3$.
\end{remark}

\section{The nonlinear remainder}
\label{sec:Remainder}

The test function $W(u)$ is far from linear in $u$ as a function on the cube. In this section we show that it is linear up to a small relative error once its component in the low Haar cross is removed and the primitive is taken. Let
\begin{equation}\label{eq:LowCross}
    \Hlo\coloneqq\operatorname{span}\bigl\{\text{tensor Haar functions of level at most }h-1\bigr\}\subset L^2(\cube),
\end{equation}
a finite-dimensional space, which we call the \emph{low Haar cross}, and let $\Plo$ be the orthogonal projection of $L^2(\cube)$ onto $\Hlo$ and $\Phig\coloneqq\mathrm{Id}-\Plo$. Since $q_R$ is a combination of tensor Haar functions of level $h+2d$, we have
\begin{equation}\label{eq:WlinHigh}
    \Plo\Wlin(u)=0\qquad\text{and}\qquad\Phig\Wlin(u)=\Wlin(u)\qquad\text{for every }u.
\end{equation}
The component in the low Haar cross has to be removed because of coincidences of depths. If, in a coordinate $k\ge2$, the largest side depth among the shell boxes in a product $\prod_Rq_R$ is attained by an even number of them, then the factors in that coordinate multiply to $\pm\one_I$, where $I$ is the common deepest side, because $q_I^2=\one_I$ and the factors of smaller depth are constant on $I$ (\cref{lem:CellMonomial}). The primitive of $q_I$ vanishes outside $I$, but the primitive of $\one_I$ equals $|I|$ to the left of $I$, so the primitive of such a product does not vanish outside its support; the part far from $I$ comes from the constant and the coarse Haar components of $\one_I$. The projection $\Phig$ removes all tensor Haar components of level below $h$. It therefore retains only components whose depths fall short of the depths of the deepest sides by at most $e+2d$ in total, where $e$ is the excess of the join of the shapes involved (\cref{sec:CellMonomials}), and \cref{lem:CellLocation} turns this into exponential decay in $e$.

\begin{theorem}[Remainder estimate]
    \label{thm:Remainder}
    There is a constant $t_d^*\in(0,\tfrac12]$ depending only on $d$ such that the following holds for every $t_d\in(0,t_d^*]$ in \eqref{eq:Parameters} and every $h\in3\N$, with constants $C_d,C_d'$ that depend only on $d$ and not on $t_d$. If $U$ is uniformly distributed on $[-1,1]^{\Shell}$, then
    \begin{equation}\label{eq:RemainderBound}
        \E_{U}\paren[\big]{\norm{\Jst\Phig W(U)-\Flin(U)}_2^2}\le C_d\,\lambda^2(h+1)^{d-1}v^2\,(\theta^4+\lambda^4)
        \le C_d'\,\theta^4\;\E_{U}\paren[\big]{\norm{\Flin(U)}_2^2}.
    \end{equation}
\end{theorem}

The second inequality follows from the first by \eqref{eq:FlinNorm}, whose constants do not depend on $t_d$, and $\lambda\le\theta$. Since $\theta^4=t_d^4(h+1)^{-2(d-2)}$, the relative root-mean-square error is $O_d((h+1)^{-(d-2)})$, which tends to zero as $h\to\infty$ because $d\ge3$.

From now on we fix $t_d\coloneqq t_d^*$, with $t_d^*$ as chosen in \eqref{eq:tdChoice}; the constants $c_d,C_d$ appearing later may depend on this choice, whereas the constants in \eqref{eq:RemainderBound} do not.

\paragraph{Strategy.}
We expand $W(u)$ into monomials in the coefficients and use their orthogonality to reduce the mean-square error to a sum of squared norms. After removal of the component in the low Haar cross, the cell-location estimate gives exponential decay in the join excess, up to polynomial factors. We then separate monomials lying in one block from those meeting at least three blocks; their total contributions are bounded by factors $O_d(\theta^4)$ and $O_d(\lambda^4)$, respectively, relative to the mean-square size of the linear part.

\subsection{Cell monomials and joins}
\label{sec:CellMonomials}

Substituting \eqref{eq:BlockFunction} and \eqref{eq:LayerFunction} into \eqref{eq:TestFunctionExpansion} and expanding all products, we write $W(u)$ as a sum of monomials in the coefficients. Each term is indexed by a set $\Gamma\subseteq\Shell$ of shell boxes with pairwise distinct shapes, one from each layer function that occurs; we write $\calT(\Gamma)\coloneqq\{\shape(R):R\in\Gamma\}$ for its set of shapes and $S(\Gamma)\coloneqq\{r_1:r\in\calT(\Gamma)\}$ for the set of blocks it meets. We call $\Gamma$ \emph{admissible} if $|S(\Gamma)|$ is odd and $|\calT(\Gamma)\cap\Bl a|$ is odd for every $a\in S(\Gamma)$, and we denote the set of admissible $\Gamma$ by $\cellmon$. For $\Gamma\in\cellmon$ we put
\begin{equation}\label{eq:CellMonomial}
    c_\Gamma\coloneqq\varepsilon^{|S(\Gamma)|}\prod_{a\in S(\Gamma)}\widehat\Phi_a\bigl(\calT(\Gamma)\cap\Bl a\bigr),
    \qquad
    m_\Gamma(u)\coloneqq\prod_{R\in\Gamma}u_R,
    \qquad
    X_\Gamma\coloneqq\prod_{R\in\Gamma}q_R .
\end{equation}
Then, by \cref{lem:WalshTanh}(i), which shows that only odd sets contribute in \eqref{eq:BlockFunction},
\begin{equation}\label{eq:WExpansion}
    W(u)=\sum_{\Gamma\in\cellmon}c_\Gamma\,m_\Gamma(u)\,X_\Gamma,
    \qquad
    \Wlin(u)=\sum_{\Gamma\in\cellmon,\ |\Gamma|=1}c_\Gamma\,m_\Gamma(u)\,X_\Gamma,
\end{equation}
where the second identity restates \cref{lem:LinearPart}; note that $|\Gamma|=1$ means $|S(\Gamma)|=1$ and $|\calT(\Gamma)|=1$, and that every other admissible $\Gamma$ has $|\Gamma|\ge3$. We call the terms $c_\Gamma m_\Gamma(u)X_\Gamma$ with $\Gamma\in\cellmon$ the \emph{cell monomials} of $W$; for $X_\Gamma\ne0$, the dyadic box $\bigcap_{R\in\Gamma}R$, on which $X_\Gamma$ is supported, is called the \emph{cell} of $\Gamma$. Each $m_\Gamma$ is squarefree, and distinct admissible sets give distinct monomials. If $U$ is uniformly distributed on $[-1,1]^{\Shell}$, then $\E(m_\Gamma(U)m_{\Gamma'}(U))=0$ for $\Gamma\ne\Gamma'$, because some coefficient occurs exactly once in the product, and $\E(m_\Gamma(U)^2)=3^{-|\Gamma|}\le1$. Consequently, for every bounded linear operator $\calO$ on $L^2(\cube)$,
\begin{equation}\label{eq:MonomialOrthogonality}
    \E_{U}\paren[\big]{\norm[\big]{\calO\bigl(W(U)-\Wlin(U)\bigr)}_2^2}
    =\sum_{\Gamma\in\cellmon,\ |\Gamma|\ge3}c_\Gamma^2\,3^{-|\Gamma|}\norm{\calO X_\Gamma}_2^2
    \le\sum_{\Gamma\in\cellmon,\ |\Gamma|\ge3}c_\Gamma^2\,\norm{\calO X_\Gamma}_2^2 .
\end{equation}
We apply this with $\calO=\Jst\Phig$, which is bounded since $\Jst$ is bounded on $L^2(\cube)$; by \eqref{eq:WlinHigh}, the left-hand side of \eqref{eq:MonomialOrthogonality} is then the left-hand side of \eqref{eq:RemainderBound}. To bound the right-hand side of \eqref{eq:MonomialOrthogonality}, we describe the functions $X_\Gamma$.

For a set of shapes $\calT\subseteq\Shapes$, the \emph{join} $j=j(\calT)\in(3\N_0)^d$ is the coordinatewise maximum $j_k\coloneqq\max_{r\in\calT}r_k$, and the \emph{excess} is $e(\calT)\coloneqq|j|_1-h$. Since $|r|_1=h$ for every shape and $r_k\le h$, we have $0\le e(\calT)\le(d-1)h$.

\begin{lemma}[Structure of cell monomials]
    \label{lem:CellMonomial}
    Let $\calT=\calT(\Gamma_0)$ for some $\Gamma_0\in\cellmon$, with join $j$ and excess $e$.
    \begin{enumerate}[label=\textup{(\alph*)}]
        \item The admissible sets $\Gamma$ with $\calT(\Gamma)=\calT$ and $X_\Gamma\ne0$ are in bijection with the dyadic boxes $I=I_1\times\dots\times I_d$ with $I_k\in\dy_{j_k}$, via $\Gamma\mapsto I\coloneqq\bigcap_{R\in\Gamma}R$. In particular there are exactly $2^{h+e}$ of them.
        \item If $\Gamma$ corresponds to $I$ in \textup{(a)}, then
        \begin{equation}\label{eq:CellFactorization}
            X_\Gamma=\pm\,q_{I_1}\otimes g_2\otimes\dots\otimes g_d,
            \qquad
            g_k=\begin{cases}q_{I_k}&\text{if }\#\{r\in\calT:r_k=j_k\}\text{ is odd},\\ \one_{I_k}&\text{otherwise.}\end{cases}
        \end{equation}
    \end{enumerate}
\end{lemma}

\begin{proof}
    (a) If $X_\Gamma\ne0$, the boxes in $\Gamma$ have a common point, so in every coordinate $k$ their sides are pairwise nested, and their intersection is a dyadic box $I$ whose $k$-th side is the deepest of these sides, of depth $\max_{r\in\calT}r_k=j_k$. Conversely, given a dyadic box $I$ with sides of depths $j_k$, every shape $r\in\calT$ satisfies $r\le j$, and the unique shell box $R_r$ of shape $r$ containing $I$ is the product of the ancestors of the sides of $I$ of depths $r_k$. The set $\Gamma\coloneqq\{R_r:r\in\calT\}$ is admissible, has $\calT(\Gamma)=\calT$ and $\bigcap_{R\in\Gamma}R=I$, and $X_\Gamma\ne0$ because $|X_\Gamma|=\prod_{R\in\Gamma}\one_R=\one_I$; the two maps are inverse to each other. The number of dyadic boxes with sides of depths $j_1,\dots,j_d$ is $2^{|j|_1}=2^{h+e}$.

    (b) Write $R_r$ for the shell box of shape $r$ in $\Gamma$ and $R_{r,k}$ for its $k$-th side, so that $X_\Gamma=\bigotimes_{k=1}^d\prod_{r\in\calT}q_{R_{r,k}}$. Fix a coordinate $k$. The product $\prod_{r\in\calT}q_{R_{r,k}}$ vanishes outside $I_k$, since it contains the factor $q_{I_k}$ (some $r$ has $r_k=j_k$). Let $x_k\in I_k$. For $r\in\calT$ with $r_k<j_k$ we have $r_k\le j_k-3$, so $I_k$ is contained in an interval of $\dy_{r_k+3}$ inside $R_{r,k}$, on which $q_{R_{r,k}}$ is constant and equal to $\pm1$. For $r\in\calT$ with $r_k=j_k$ we have $R_{r,k}=I_k$ and the factor is $q_{I_k}$. Since $q_{I_k}^2=\one_{I_k}$, the product $\prod_{r\in\calT}q_{R_{r,k}}(x_k)$ equals $\pm q_{I_k}(x_k)$ or $\pm\one_{I_k}(x_k)$ according to the parity of $\#\{r:r_k=j_k\}$. In the first coordinate, $j_1=\max S(\Gamma)$ is the block of largest index met by $\Gamma$, and the shapes with $r_1=j_1$ are exactly the elements of $\calT\cap\Bl{j_1}$, whose number is odd by admissibility.
\end{proof}

\subsection{The cell-location estimate}
\label{sec:CellLocation}

The following lemma bounds the total contribution of all monomials with a fixed set of shapes.

\begin{lemma}[Cell-location estimate]
    \label{lem:CellLocation}
    Let $\calT$ be as in \cref{lem:CellMonomial}, with excess $e$. Then
    \begin{equation*}
        \sum_{\Gamma\in\cellmon:\ \calT(\Gamma)=\calT}\norm{\Jst\Phig X_\Gamma}_2^2\le C_d\,v^2\,2^{-e}\,(e+1)^{2(d-1)}.
    \end{equation*}
\end{lemma}

\begin{proof}
    Sets $\Gamma$ with $X_\Gamma=0$ contribute nothing. Let $\Gamma$ with $X_\Gamma\ne0$ correspond to the box $I$ in \cref{lem:CellMonomial}(a), and write $X_\Gamma$ as in \eqref{eq:CellFactorization}. Call a coordinate $k\ge2$ \emph{modulated} if $g_k=q_{I_k}$ and \emph{plain} otherwise. We expand every factor in the Haar basis. The factor $q_{I_1}$ and every modulated factor $q_{I_k}$ are combinations of four Haar functions of level $j_1+2$ and $j_k+2$, respectively. A plain factor $\one_{I_k}$ has the expansion
    \begin{equation}\label{eq:IndicatorExpansion}
        \one_{I_k}=2^{-j_k}\one+\sum_{b=0}^{j_k-1}\sigma_b\,2^{b-j_k}\,h_{J_b},
    \end{equation}
    where $J_b$ is the ancestor of $I_k$ of depth $b$ and $\sigma_b\in\{-1,1\}$ records on which half of $J_b$ the interval $I_k$ lies. For a \emph{component vector} $\beta=(\beta_2,\dots,\beta_d)$ with $\beta_k=j_k$ for modulated $k$ and $\beta_k\in\{-1,0,\dots,j_k-1\}$ for plain $k$, let $g_k^{(\beta_k)}$ be the corresponding term: $g_k^{(j_k)}\coloneqq q_{I_k}$ for modulated $k$, and $g_k^{(-1)}\coloneqq2^{-j_k}\one$, $g_k^{(b)}\coloneqq\sigma_b2^{b-j_k}h_{J_b}$ for plain $k$. Then
    \begin{equation*}
        X_\Gamma=\pm\sum_\beta X_\Gamma^{(\beta)},\qquad X_\Gamma^{(\beta)}\coloneqq q_{I_1}\otimes g_2^{(\beta_2)}\otimes\dots\otimes g_d^{(\beta_d)},
    \end{equation*}
    and $X_\Gamma^{(\beta)}$ is a combination of tensor Haar functions of the common level
    \begin{equation*}
        \lev(\beta)=(j_1+2)+\sum_{k\ge2}\beta_k+2\,\#\{k\ge2:k\text{ modulated}\}.
    \end{equation*}
    Hence $\Phig X_\Gamma=\pm\sum_{\beta:\ \lev(\beta)\ge h}X_\Gamma^{(\beta)}$, and every retained component vector satisfies
    \begin{equation}\label{eq:Retention}
        j_1+\sum_{k\ge2}\beta_k\ \ge\ h-2d .
    \end{equation}
    The deficits $\delta_k\coloneqq j_k-\beta_k\ge0$ of a retained component vector satisfy $\sum_{k\ge2}\delta_k=(h+e-j_1)-\sum_{k\ge2}\beta_k\le e+2d$ by \eqref{eq:Retention}. The number $N_\calT$ of retained component vectors is therefore at most the number of vectors in $\N_0^{d-1}$ with coordinate sum at most $e+2d$, so that
    \begin{equation}\label{eq:ComponentCount}
        N_\calT\le\binom{e+3d-1}{d-1}\le C_d(e+1)^{d-1}.
    \end{equation}

    Next we bound the primitives. By \eqref{eq:psiProperties} and \eqref{eq:HaarPrimitives}, $\norm{\Jst q_{I_1}}_2^2=2^{-3j_1}/192$, and for $k\ge2$
    \begin{equation*}
        \norm[\big]{\Jst g_k^{(j_k)}}_2^2=\frac{2^{-3j_k}}{192},\qquad
        \norm[\big]{\Jst g_k^{(-1)}}_2^2=\frac{2^{-2j_k}}{3},\qquad
        \norm[\big]{\Jst g_k^{(b)}}_2^2=2^{2b-2j_k}\cdot\frac{2^{-3b}}{12},\qquad0\le b<j_k,
    \end{equation*}
    so that in all three cases $\norm{\Jst g_k^{(\beta_k)}}_2^2\le2^{-2j_k-\beta_k}$. Since $\Jst$ acts coordinatewise on tensor products,
    \begin{equation*}
        \norm[\big]{\Jst X_\Gamma^{(\beta)}}_2^2\le2^{-3j_1-\sum_{k\ge2}(2j_k+\beta_k)}=2^{-2(h+e)-j_1-\sum_{k\ge2}\beta_k}\le2^{2d}\,2^{-3h-2e}=2^{2d}\,v^3\,2^{-2e},
    \end{equation*}
    where we used $\sum_{k\ge2}2j_k=2(h+e)-2j_1$ and \eqref{eq:Retention}. This bound does not depend on the box $I$. The components $X_\Gamma^{(\beta)}$ for different $\beta$ need not have orthogonal primitives, so we use the Cauchy--Schwarz inequality:
    \begin{equation*}
        \norm{\Jst\Phig X_\Gamma}_2^2\le N_\calT\sum_{\beta:\ \lev(\beta)\ge h}\norm[\big]{\Jst X_\Gamma^{(\beta)}}_2^2\le N_\calT^2\,2^{2d}\,v^3\,2^{-2e}.
    \end{equation*}
    Summing over the $2^{h+e}=v^{-1}2^{e}$ admissible sets $\Gamma$ with $\calT(\Gamma)=\calT$ and $X_\Gamma\ne0$ (the others contribute nothing) and using \eqref{eq:ComponentCount} gives the claim.
\end{proof}

\subsection{Proof of \texorpdfstring{\cref{thm:Remainder}}{the remainder estimate}}
\label{sec:RemainderProof}

We first record two counting bounds.

\begin{lemma}[Counting shapes and joins]
    \label{lem:Counting}
    Let $j\in(3\N_0)^d$ with $|j|_1=h+e$, where $0\le e\le(d-1)h$.
    \begin{enumerate}[label=\textup{(\alph*)}]
        \item The number of shapes $r\in\Shapes$ with $r\le j$ is at most $C_d(e+1)^{d-1}$, and the number of those with $r_1=j_1$ is at most $C_d(e+1)^{d-2}$.
        \item The number of vectors $j'\in(3\N_0)^d$ with $|j'|_1=h+e$ is at most $C_d(h+1)^{d-1}$.
    \end{enumerate}
\end{lemma}

\begin{proof}
    (a) A shape $r\le j$ is determined by the vector $j-r\in(3\N_0)^d$, whose coordinate sum is $e$. The number of such vectors is at most $\binom{e/3+d-1}{d-1}\le C_d(e+1)^{d-1}$, and at most $\binom{e/3+d-2}{d-2}\le C_d(e+1)^{d-2}$ if the first coordinate is prescribed to be zero.

    (b) The number of such $j'$ is $\binom{(h+e)/3+d-1}{d-1}\le C_d(h+e+1)^{d-1}\le C_d(h+1)^{d-1}$, since $e\le(d-1)h$.
\end{proof}

\begin{proof}[Proof of \cref{thm:Remainder}]
By \eqref{eq:MonomialOrthogonality} with $\calO=\Jst\Phig$, the identity \eqref{eq:WlinHigh}, and the fact that $c_\Gamma$ depends only on $\calT(\Gamma)$,
\begin{equation}\label{eq:RemainderGrouped}
\begin{aligned}
    \E_{U}\paren[\big]{\norm{\Jst\Phig W(U)-\Flin(U)}_2^2}
    &\le\sum_{\calT}c_\calT^2\sum_{\Gamma\in\cellmon:\ \calT(\Gamma)=\calT}\norm{\Jst\Phig X_\Gamma}_2^2\\
    &\le\sum_{\calT}c_\calT^2\;C_dv^2\,2^{-e(\calT)}\,(e(\calT)+1)^{2(d-1)},
\end{aligned}
\end{equation}
where $\calT$ ranges over the sets of shapes of admissible sets with at least three elements, $c_\calT$ denotes the common value of $c_\Gamma$ for $\calT(\Gamma)=\calT$, and we used \cref{lem:CellLocation}. We group the sets $\calT$ according to their join $j$ and distinguish two classes. Since $|S(\Gamma)|$ is odd for admissible $\Gamma$, either $\calT$ meets a single block, and then $|\calT|\ge3$ because $|\Gamma|\ge3$ (class~I), or it meets at least three blocks (class~II).

\paragraph{Class I: one block.}
Here $\calT\subseteq\Bl a$ for a single block $a$, with $|\calT|\ge3$ odd, and $c_\calT=\varepsilon\widehat\Phi_a(\calT)$. The join $j$ of $\calT$ satisfies $j_1=a$, and $\calT$ is a subset of $\calI_j\coloneqq\{r\in\Bl{j_1}:r\le j\}$, whose cardinality is at most $C_d(e+1)^{d-2}$ by \cref{lem:Counting}(a). By the second inequality in \eqref{eq:WalshActivity}, applied to $F=\Phi_{j_1}$ and $\calE=\calI_j$,
\begin{equation*}
    \sum_{\calT\subseteq\calI_j,\ |\calT|\ge3}c_\calT^2\le\varepsilon^2\theta^6|\calI_j|^3\le C_d\,\varepsilon^2\theta^6(e+1)^{3(d-2)} .
\end{equation*}
For each excess $e$ there are at most $C_d(h+1)^{d-1}$ joins by \cref{lem:Counting}(b), and every excess satisfies $0\le e\le(d-1)h$. Hence the contribution of class I to the right-hand side of \eqref{eq:RemainderGrouped} is at most
\begin{equation}\label{eq:ClassI}
\begin{aligned}
    \MoveEqLeft\sum_{e=0}^{(d-1)h}C_d(h+1)^{d-1}\cdot\varepsilon^2\theta^6(e+1)^{3(d-2)}\cdot v^22^{-e}(e+1)^{2(d-1)}\\
    &\le C_d\,\varepsilon^2\theta^6(h+1)^{d-1}v^2\sum_{e\ge0}2^{-e}(e+1)^{5d-8}
    \le C_d\,\lambda^2\theta^4(h+1)^{d-1}v^2 .
\end{aligned}
\end{equation}

\paragraph{Class II: at least three blocks.}
Here $|S|\ge3$, where $S$ is the set of blocks met by $\calT$, and $c_\calT^2=\prod_{a\in S}\varepsilon^2\widehat\Phi_a(\calT\cap\Bl a)^2$. Fix a join $j$ with excess $e$ and put $\calI'_j\coloneqq\{r\in\Shapes:r\le j\}$, so that $|\calI'_j|\le C_d(e+1)^{d-1}$ by \cref{lem:Counting}(a). For a block $a$ let
\begin{equation*}
    \alpha_a(j)\coloneqq\varepsilon^2\sum_{\emptyset\ne T\subseteq\Bl a\cap\calI'_j}\widehat\Phi_a(T)^2\le\varepsilon^2\theta^2\,|\Bl a\cap\calI'_j|,
    \qquad
    \alpha(j)\coloneqq\sum_a\alpha_a(j)\le\varepsilon^2\theta^2|\calI'_j|\le C_d\,\lambda^2(e+1)^{d-1},
\end{equation*}
where we used the first inequality in \eqref{eq:WalshActivity}. Every class~II set $\calT$ with join $j$ has the form $\bigcup_{a\in S}T_a$ with nonempty $T_a\subseteq\Bl a\cap\calI'_j$. Summing $c_\calT^2$ over all such choices with a fixed $S$ gives at most $\prod_{a\in S}\alpha_a(j)$, and summing over all $S$ with $|S|\ge3$ gives at most
\begin{equation*}
    \sum_{k\ge3}\ \sum_{|S|=k}\prod_{a\in S}\alpha_a(j)\le\sum_{k\ge3}\frac{\alpha(j)^k}{k!}\le\frac{\alpha(j)^3}{6}\mathrm{e}^{\alpha(j)},
\end{equation*}
because the $k$-th elementary symmetric polynomial of nonnegative numbers is at most the $k$-th power of their sum divided by $k!$. Using \cref{lem:Counting}(b) once more and recalling that the excess of every set of shapes satisfies $0\le e\le(d-1)h$, the contribution of class II to \eqref{eq:RemainderGrouped} is at most
\begin{equation}\label{eq:ClassII}
    C_d(h+1)^{d-1}v^2\sum_{e=0}^{(d-1)h}2^{-e}(e+1)^{2(d-1)}\cdot\lambda^6(e+1)^{3(d-1)}\exp\bigl(C_d\lambda^2(e+1)^{d-1}\bigr) .
\end{equation}
For every $e$ in this range, $e+1\le(d-1)h+1\le d(h+1)$, so the definition \eqref{eq:Parameters} of $\theta$ and $\lambda$ gives
\begin{equation*}
    \lambda^2(e+1)^{d-1}=\varepsilon^2t_d^2\,(h+1)^{-(d-2)}(e+1)^{d-2}\,(e+1)\le\varepsilon^2t_d^2\,d^{\,d-2}\,(e+1).
\end{equation*}
We now choose
\begin{equation}\label{eq:tdChoice}
    t_d^*\in(0,\tfrac12]\qquad\text{with}\qquad C_d\,\varepsilon^2(t_d^*)^2d^{\,d-2}\le\tfrac14\log2,
\end{equation}
where $C_d$ is the constant in the exponent in \eqref{eq:ClassII}; it comes from \cref{lem:Counting}(a) and does not depend on $t_d$. Then, for every $t_d\in(0,t_d^*]$, the exponential factor is at most $2^{(e+1)/4}$ for all $0\le e\le(d-1)h$, so the sum in \eqref{eq:ClassII}, without the factor $\lambda^6$, is at most ${\sum_{e\ge0}2^{-e}(e+1)^{5(d-1)}2^{(e+1)/4}}$, a convergent series whose value depends only on $d$; hence the contribution of class II is at most $C_d\lambda^6(h+1)^{d-1}v^2$. Together with \eqref{eq:ClassI} this proves the first inequality in \eqref{eq:RemainderBound}.
\end{proof}

\section{Sampling and the exact linear response}
\label{sec:Sampling}

We choose the point set and the matrix $B$ that transforms a standard Gaussian vector $z\in\R^M$ into the weight vector $\eta=Bz\in\R^n$. Three properties are needed: the linear part of the pairing of the weighted counting function $G_z$ with $W(u)$ is exactly diagonal in $u$; the pairing does not see the low Haar cross, so that \cref{thm:Remainder} controls the nonlinear remainder; and $BB^*\preceq I_n$, which gives the variance bound required in \cref{lem:GaussianCover}. The first two properties rest on a weighted point set for which discrete and continuous $L^2$-norms are equivalent on a space containing all $T_R$ and the primitives of the low Haar cross (\cref{sec:AffineSpace}), and which also controls the mean-square remainder at the points (\cref{sec:SamplingLemma}); \cref{sec:Response} then defines $B$ by a projection and an inverse Gram matrix.

\subsection{A space of piecewise multiaffine functions}
\label{sec:AffineSpace}

For $j\in\N_0$ let $V_j$ be the space of functions on $[0,1)$ that are affine on every interval of $\dy_j$; thus $\dim V_j=2^{j+1}$ and $V_0\subset V_1\subset\cdots$. Put $D_0\coloneqq V_0$ and $D_j\coloneqq V_j\ominus V_{j-1}$ (the orthogonal complement of $V_{j-1}$ in $V_j$) for $j\ge1$, so that $\dim D_j=2^j$ for $j\ge1$ and the spaces $D_j$ are mutually orthogonal in $L^2[0,1)$. For $\ell\in\N_0^d$ let $D_\ell\coloneqq D_{\ell_1}\otimes\dots\otimes D_{\ell_d}\subset L^2([0,1)^d)$, and define
\begin{equation}\label{eq:VH}
    H\coloneqq h+3d,\qquad
    \VH\coloneqq\bigoplus_{\ell\in\N_0^d,\ |\ell|_1\le H}D_\ell .
\end{equation}
Every element of $\VH$ has a unique representative that is multiaffine, that is, affine in each coordinate separately, on each of the $2^{dH}$ half-open dyadic boxes with sides in $\dy_H$; we always use this representative, which is defined at every point of $[0,1)^d$.

For a finite-dimensional space $X$ of pointwise defined functions on a set $\Omega$ with a probability measure $\mu$, the \emph{inverse Christoffel function} is $L_X(p)\coloneqq\sum_k|\phi_k(p)|^2$, where $(\phi_k)$ is an orthonormal basis of $X$ in $L^2(\mu)$. It does not depend on the choice of the basis, since $L_X(p)$ is the squared norm of the evaluation functional $f\mapsto f(p)$ on $X$. It is nondecreasing in $X$, additive for orthogonal direct sums, and multiplicative for tensor products of spaces on product measure spaces. \Cref{thm:Sampling} requires a uniform bound on the inverse Christoffel function of the sampling space; in the proof of \cref{lem:SamplingLemma}, the lower bound in part~(c) of the following lemma keeps the weights bounded by $C_d/n$ after a change of density, while part~(b) controls the number of points.

\begin{lemma}[Inverse Christoffel function of $\VH$]
    \label{lem:AffineSpace}
    \begin{enumerate}[label=\textup{(\alph*)}]
        \item There are absolute constants $c_0,C_0>0$ such that $c_0\,2^j\le L_{D_j}(t)\le C_0\,2^j$ for all $j\in\N_0$ and $t\in[0,1)$.
        \item $c_d\,2^H(H+1)^{d-1}\le\dim\VH\le C_d\,2^H(H+1)^{d-1}$; in particular $c_dM\le\dim\VH\le C_dM$.
        \item $c_d\dim\VH\le L_{\VH}(p)\le C_d\dim\VH$ for all $p\in[0,1)^d$.
    \end{enumerate}
\end{lemma}

\begin{proof}
    (a) For $I=[x_I,x_I+2^{-j})\in\dy_j$, the functions $2^{j/2}\one_I$ and $2^{j/2}\sqrt3\,(2^{j+1}(t-x_I)-1)\one_I$ form an orthonormal basis of the functions that are affine on $I$ and vanish outside $I$, so for $t\in I$,
    \begin{equation*}
        L_{V_j}(t)=2^j\bigl(1+12(s-\tfrac12)^2\bigr),\qquad s\coloneqq2^j(t-x_I)\in[0,1).
    \end{equation*}
    For $j=0$ this gives $1\le L_{D_0}=L_{V_0}\le4$. For $j\ge1$ and $t$ in the left half $I$ of $I'\in\dy_{j-1}$ we have $2^{j-1}(t-x_{I'})=s/2$, hence
    \begin{equation*}
        L_{D_j}(t)=L_{V_j}(t)-L_{V_{j-1}}(t)=2^{j-1}\bigl(4+21s^2-18s\bigr)\in\bigl[\tfrac17\,2^{j-1},\,7\cdot2^{j-1}\bigr],
    \end{equation*}
    and the case of the right half is symmetric.

    (b) Since $2^{|\ell|_1}\le\dim D_\ell=\prod_k\dim D_{\ell_k}\le2^d2^{|\ell|_1}$, and since $\sum_{|\ell|_1\le H}2^{|\ell|_1}=\sum_{b=0}^H2^b\binom{b+d-1}{d-1}$ lies between $2^H\binom{H+d-1}{d-1}$ and $2^{H+1}\binom{H+d-1}{d-1}$, we get $\dim\VH\asymp2^H(H+1)^{d-1}$. By \eqref{eq:VH}, \eqref{eq:ShapeCount} and \eqref{eq:MDefinition}, this is $\asymp(h+1)^{d-1}/v\asymp M$.

    (c) By additivity and multiplicativity, $L_{\VH}=\sum_{|\ell|_1\le H}\prod_kL_{D_{\ell_k}}$, and (a) gives $c_0^d\sum_{|\ell|_1\le H}2^{|\ell|_1}\le L_{\VH}\le C_0^d\sum_{|\ell|_1\le H}2^{|\ell|_1}$, which is $\asymp\dim\VH$ by the proof of (b).
\end{proof}

Let $\Low\coloneqq\Jst(\Hlo)$ be the space of primitives of the low Haar cross \eqref{eq:LowCross}; since $\Jst$ is injective on $L^2(\cube)$, $\dim\Low=\dim\Hlo$.

\begin{lemma}[Space inclusions]
    \label{lem:SpaceInclusions}
    \begin{enumerate}[label=\textup{(\alph*)}]
        \item $T_R\in\VH$ for every $R\in\Shell$.
        \item $\Low\subset\VH$.
        \item $T_R\perp\Low$ for every $R\in\Shell$.
    \end{enumerate}
\end{lemma}

\begin{proof}
    (a) By \cref{lem:qProperties}, $\psi_{R_k}$ is affine on every interval of $\dy_{r_k+3}$, so $\psi_{R_k}\in V_{r_k+3}=\bigoplus_{\ell\le r_k+3}D_\ell$ and $T_R\in\bigoplus_{\ell_k\le r_k+3\ \forall k}D_\ell\subset\VH$, because $\sum_k(r_k+3)=h+3d=H$.

    (b) Let $\varphi=\varphi_1\otimes\dots\otimes\varphi_d$ be a tensor Haar function of level at most $h-1$. By \eqref{eq:HaarPrimitives}, $\Jst\one=1-t\in V_0$ and $\Jst h_J=-\tau_J\in V_{\depth(J)+1}$, so $\Jst\varphi_k\in V_{\lev(\varphi_k)+1}$ and $\Jst\varphi\in\bigoplus_{\ell_k\le\lev(\varphi_k)+1\ \forall k}D_\ell$, where $\sum_k(\lev(\varphi_k)+1)\le h-1+d\le H$. Since $\Low$ is spanned by these $\Jst\varphi$, (b) follows.

    (c) With $\varphi$ as in (b), $\inp{T_R}{\Jst\varphi}=v^{-1}\prod_k\inp{\psi_{R_k}}{\Jst\varphi_k}$. If some $\varphi_k=\one$, then $\Jst\varphi_k=1-t$ is affine on $[0,1)$ and the $k$-th factor vanishes by \cref{lem:qProperties}. Otherwise $\varphi_k=h_{J_k}$ for every $k$, and $\sum_k\depth(J_k)\le h-1<h=\sum_kr_k$, so $\depth(J_k)\le r_k-1$ for some $k$. Then $\tau_{J_k}$ is affine on each half of $J_k$, which are intervals of depth at most $r_k$, and $R_k$ is either disjoint from $J_k$ or contained in one of its halves; in both cases $\inp{\psi_{R_k}}{\tau_{J_k}}=0$ by \cref{lem:qProperties}, and (c) follows by linearity.
\end{proof}

\subsection{Sampling with one additional function}
\label{sec:SamplingLemma}

For $u\in[-1,1]^{\Shell}$ let
\begin{equation}\label{eq:RemainderFunction}
    \Delta_u\coloneqq\Jst\Phig W(u)-\Flin(u),\qquad
    w(p)\coloneqq\E_{U}\paren[\big]{|\Delta_U(p)|^2},\qquad p\in\cube,
\end{equation}
with $U$ uniform on $[-1,1]^{\Shell}$. By \eqref{eq:WExpansion} and the orthogonality of the monomials $m_\Gamma(U)$,
\begin{equation*}
    w(p)=\sum_{\Gamma\in\cellmon,\ |\Gamma|\ge3}c_\Gamma^2\,3^{-|\Gamma|}\,\bigl|(\Jst\Phig X_\Gamma)(p)\bigr|^2
\end{equation*}
is a finite sum; each $\Phig X_\Gamma=X_\Gamma-\Plo X_\Gamma$ is bounded, so $\Jst\Phig X_\Gamma$ is Lipschitz continuous, and therefore $w$ is continuous on $\cube$. Integrating this identity over the cube and using \eqref{eq:MonomialOrthogonality} with $\calO=\Jst\Phig$ and \cref{thm:Remainder} (with $\lambda\le\theta$) gives
\begin{equation}\label{eq:ScalarIntegral}
    \int_{\cube}w(p)\,dp=\E_{U}\paren[\big]{\norm{\Delta_U}_2^2}\le C_d\,\lambda^2(h+1)^{d-1}v^2\,\theta^4 .
\end{equation}
To transfer this mean-square remainder bound to the sample points, we add the single function $\sqrt w$ to the sampling space.

Let $\Omega\subset(0,1)^d$ be the set of points of $[0,1)^d$ none of whose coordinates is an integer multiple of $2^{-H}$. It has full measure and is the union of the open dyadic boxes with sides of depth $H$; every function in $\VH$ is continuous on each of these open boxes, which is used to make the sample points distinct.

\begin{lemma}[Sampling lemma]
    \label{lem:SamplingLemma}
    Let $w\colon\cube\to[0,\infty)$ be continuous. There are $n$ distinct points $P_1,\dots,P_n\in\Omega$ and weights $a_1,\dots,a_n>0$ such that
    \begin{equation}\label{eq:SamplingSize}
        c_dM\le n\le C_dM,\qquad a_i\le C_d/n,\qquad i\le n,
    \end{equation}
    and, with absolute constants $c_\star,C_\star>0$,
    \begin{align}
        c_\star\norm{f}_2^2&\le\sum_{i=1}^na_i|f(P_i)|^2\le C_\star\norm{f}_2^2\qquad\text{for all }f\in\VH,\label{eq:SamplingNorm}\\
        \sum_{i=1}^na_i\,w(P_i)&\le C_\star\int_{\cube}w(p)\,dp.\label{eq:SamplingScalar}
    \end{align}
\end{lemma}

\begin{proof}
    We use the change-of-density argument from the proof of \cite[Theorem~6.4]{DPSTT21} (Theorem~6.3 in arXiv:2001.10636v1); see also \cite[Remark~3.1]{LimonovaTemlyakov22}. We include the details to obtain distinct points and the required bound on the weights. Let $\calU$ be the span of $\VH$ and $\sqrt w$, considered as functions on $\Omega$. All functions in $\calU$ are continuous on each open dyadic box of depth $H$, so a function in $\calU$ that vanishes almost everywhere vanishes on $\Omega$, and $L^2(\Omega)$ with Lebesgue measure induces an inner product on $\calU$; note that $\dim\calU\in\{\dim\VH,\dim\VH+1\}$. Let $(\phi_k)_{k\le\dim\calU}$ be an orthonormal basis of $\calU$ in $L^2(\Omega)$, and put $\varpi\coloneqq L_\calU/\dim\calU=(\dim\calU)^{-1}\sum_k|\phi_k|^2$. Then $\int_\Omega\varpi=1$, and since $L_X$ is nondecreasing in $X$, \cref{lem:AffineSpace}(c) gives $\varpi\ge L_{\VH}/\dim\calU\ge c_d\dim\VH/(\dim\VH+1)\ge c_d'>0$ on $\Omega$. The functions $u_k\coloneqq\phi_k/\sqrt\varpi$ are orthonormal in $L^2(\Omega,\varpi\,dp)$ and satisfy $\sum_k|u_k|^2=\dim\calU$ on $\Omega$. By \cref{thm:Sampling}, applied with $\mu=\varpi\,dp$, there are points $\xi^1,\dots,\xi^m\in\Omega$ with $\dim\calU\le m\le C_1\dim\calU$ such that for every $g=f/\sqrt\varpi$ with $f\in\calU$,
    \begin{equation*}
        C_2\norm{f}_2^2=C_2\norm{g}_{L^2(\mu)}^2\le\frac1m\sum_{j=1}^m\frac{|f(\xi^j)|^2}{\varpi(\xi^j)}\le C_3\norm{g}_{L^2(\mu)}^2=C_3\norm{f}_2^2 .
    \end{equation*}
    With $\tilde a_j\coloneqq1/(m\varpi(\xi^j))\le1/(mc_d')$ this is a two-sided inequality $C_2\norm{f}_2^2\le\sum_j\tilde a_j|f(\xi^j)|^2\le C_3\norm{f}_2^2$ for all $f\in\calU$, which contains \eqref{eq:SamplingNorm} and, for $f=\sqrt w$, \eqref{eq:SamplingScalar}.

    The points $\xi^j$ need not be distinct. Let $\delta_0>0$, to be chosen below. Like the functions in $\calU$, the functions $\varpi$ and $1/(m\varpi)$ are continuous on each open dyadic box of depth $H$. If several points coincide, we replace them by points in the same open box at distance less than $\delta_0$ from the original ones, chosen so that all $m$ points become distinct, and recompute the weights $\tilde a_j=1/(m\varpi(\xi^j))\le1/(mc_d')$. The functional $f\mapsto\sum_j\tilde a_j|f(\xi^j)|^2$ depends continuously on the points, and the unit sphere of $\calU$ is compact, so for $\delta_0$ small enough the two-sided inequality persists with the constants $C_2/2$ and $2C_3$. We take $n\coloneqq m$, $P_i\coloneqq\xi^i$ and $a_i\coloneqq\tilde a_i$; then \eqref{eq:SamplingNorm} and \eqref{eq:SamplingScalar} hold with $c_\star\coloneqq C_2/2$ and $C_\star\coloneqq2C_3$, and $n$ lies between $\dim\calU\ge\dim\VH\ge c_dM$ and $C_1(\dim\VH+1)\le C_dM$ by \cref{lem:AffineSpace}(b), so that $a_i\le C_d/n$.
\end{proof}

\paragraph{The point set.}
From now on, $P_1,\dots,P_n\in\Omega$ and $a_1,\dots,a_n$ denote the points and weights provided by \cref{lem:SamplingLemma} for the function $w$ in \eqref{eq:RemainderFunction}. They are deterministic: they depend on $h$ and $d$ only, and neither on $z$ nor on any realization of $U$. \Cref{sec:Conclusion} applies the hereditary-discrepancy argument to this point set.

\subsection{The Gaussian weights and the response}
\label{sec:Response}

We use the weight matrix $E\coloneqq\diag(\sqrt{a_1},\dots,\sqrt{a_n})\in\R^{n\times n}$ and the evaluation matrix $T\coloneqq\bigl(T_R(P_i)\bigr)_{i\le n,\,R\in\Shell}\in\R^{n\times M}$ of the primitives.

Let $\Lambda\in\R^{n\times\dim\Low}$ be the evaluation matrix at the points of a fixed orthonormal basis of $\Low$, so that the vectors $\Lambda b$, $b\in\R^{\dim\Low}$, are exactly the vectors $(g(P_i))_{i\le n}$ with $g\in\Low$, and let $\Pi\in\R^{n\times n}$ be the orthogonal projection onto the orthogonal complement of $\range(E\Lambda)$. The \emph{Gram matrix} is
\begin{equation}\label{eq:GramMatrix}
    \Gram\coloneqq T^*E\Pi ET\in\R^{M\times M}.
\end{equation}

\begin{lemma}[Gram matrix]
    \label{lem:Gram}
    Let $c_\star,C_\star$ be the constants of \eqref{eq:SamplingNorm}. Then $c_\star\gamma_dv\,I_M\preceq\Gram\preceq C_\star\gamma_dv\,I_M$.
\end{lemma}

\begin{proof}
    For $y\in\R^M$ put $F_y\coloneqq\sum_Ry_RT_R$, so that $ETy=(\sqrt{a_i}F_y(P_i))_i$. Since $\Pi$ projects onto $\range(E\Lambda)^\perp$,
    \begin{equation*}
        y^*\Gram y=\norm{\Pi ETy}_2^2=\min_{b}\norm{ETy-E\Lambda b}_2^2=\min_{g\in\Low}\sum_ia_i\bigl|F_y(P_i)-g(P_i)\bigr|^2 .
    \end{equation*}
    By \cref{lem:SpaceInclusions}, $F_y-g\in\VH$ and $F_y\perp g$ for every $g\in\Low$, so \eqref{eq:SamplingNorm} and \cref{lem:TOrthogonality} give $\sum_ia_i|F_y(P_i)-g(P_i)|^2\ge c_\star\norm{F_y-g}_2^2\ge c_\star\norm{F_y}_2^2=c_\star\gamma_dv\norm{y}_2^2$ for every $g$, and hence $y^*\Gram y\ge c_\star\gamma_dv\norm{y}_2^2$. Omitting the projection, $y^*\Gram y\le\norm{ETy}_2^2=\sum_ia_i|F_y(P_i)|^2\le C_\star\norm{F_y}_2^2=C_\star\gamma_dv\norm{y}_2^2$.
\end{proof}

We now define the scale $\alpha$ and the matrix $B$ by
\begin{equation}\label{eq:BDefinition}
    \alpha\coloneqq\omega_d\sqrt{n/v},\qquad
    B\coloneqq\alpha\gamma_dv\;E\Pi ET\,\Gram^{-1}\in\R^{n\times M},
\end{equation}
where $\omega_d>0$ is a constant depending only on $d$, chosen in the next lemma. The two factors of $B$ serve different purposes. The projection $\Pi$ removes the evaluated primitives of the low Haar cross: the weight vectors $Bz$ are orthogonal to every vector $(g(P_i))_{i\le n}$ with $g\in\Low$ (\cref{lem:BProperties}(a)), so that the weighted counting function does not see the low Haar cross. The inverse Gram matrix makes the pairing with the primitives $T_R$ exactly diagonal, $T^*B=\alpha\gamma_dv\,I_M$ (\cref{lem:BProperties}(b)); hence the linear part of the response pairs each coordinate $z_R$ only with $u_R$ (\cref{prop:Response}(d)), which matches \cref{sec:Entropy}, where the correlations are controlled coordinate by coordinate (\cref{prop:Correlation}).

\begin{lemma}[Properties of $B$]
    \label{lem:BProperties}
    \begin{enumerate}[label=\textup{(\alph*)}]
        \item $\Lambda^*B=0$.
        \item $T^*B=\alpha\gamma_dv\,I_M$.
        \item If $\omega_d$ is sufficiently small, depending only on $d$, then $BB^*\preceq I_n$.
    \end{enumerate}
\end{lemma}

\begin{proof}
    (a) Since $\Pi$ is symmetric and annihilates $\range(E\Lambda)$, we have $\Lambda^*E\Pi=(\Pi E\Lambda)^*=0$, and hence $\Lambda^*B=0$. (b) $T^*B=\alpha\gamma_dv\,T^*E\Pi ET\,\Gram^{-1}=\alpha\gamma_dv\,\Gram\Gram^{-1}$. (c) We have $\norm{B}_{\mathrm{op}}\le\alpha\gamma_dv\,\norm{E}_{\mathrm{op}}\norm{ET}_{\mathrm{op}}\norm{\Gram^{-1}}_{\mathrm{op}}$, with $\norm{E}_{\mathrm{op}}^2=\max_ia_i\le C_d/n$, $\norm{ET}_{\mathrm{op}}^2=\norm{T^*E^2T}_{\mathrm{op}}\le C_\star\gamma_dv$ by \eqref{eq:SamplingNorm} and \cref{lem:TOrthogonality}, and $\norm{\Gram^{-1}}_{\mathrm{op}}\le(c_\star\gamma_dv)^{-1}$ by \cref{lem:Gram}. Hence $\norm{B}_{\mathrm{op}}^2\le C_d\,\alpha^2\gamma_dv/n=C_d\gamma_d\omega_d^2$, which is at most one for small $\omega_d$; this is equivalent to $BB^*\preceq I_n$.
\end{proof}

We fix $\omega_d$ as in \cref{lem:BProperties}(c). For $z\in\R^M$, the \emph{weight vector} is $\eta\coloneqq Bz\in\R^n$, and the associated \emph{weighted counting function} is the function $G_\eta$ of \eqref{eq:DualityIntro},
\begin{equation}\label{eq:IncidenceField}
    G_z(x)\coloneqq\sum_{i=1}^n(Bz)_i\,\one_{\{P_i\le x\}},\qquad x\in\cube .
\end{equation}
For $u\in[-1,1]^{\Shell}$, the \emph{response} of the point set to the test function $W(u)$ is
\begin{equation}\label{eq:Response}
    V(u)\coloneqq B^*\bigl((\Jst W(u))(P_i)\bigr)_{i\le n}\in\R^M .
\end{equation}
The linear part of the response is exactly diagonal; the nonlinear part contributes the error $e(u)$, estimated in part~(e) of the following proposition.

\begin{proposition}[Exact linear response]
    \label{prop:Response}
    Let $s_R\coloneqq\alpha\gamma_dv^2\lambda_R$ for $R\in\Shell$, $s\coloneqq\alpha\gamma_dv^2\lambda$, and $s_*\coloneqq\min_Rs_R$, so that $s/2\le s_*\le s_R\le s$ by \cref{lem:LinearPart}. Then the following hold.
    \begin{enumerate}[label=\textup{(\alph*)}]
        \item $\inp{G_z}{W(u)}=z\cdot V(u)$ for all $z\in\R^M$ and $u\in[-1,1]^{\Shell}$.
        \item $\Plo G_z=0$ for every $z$, and $B^*\bigl((\Jst f)(P_i)\bigr)_i=B^*\bigl((\Jst\Phig f)(P_i)\bigr)_i$ for every $f\in L^2(\cube)$.
        \item $\inp{G_z}{q_R}=\alpha\gamma_dv^2\,z_R$ for all $z\in\R^M$ and $R\in\Shell$.
        \item $V(u)=\diag(s_R)\,u+e(u)$, where $e(u)\coloneqq B^*\bigl(\Delta_u(P_i)\bigr)_{i\le n}$ with $\Delta_u$ as in \eqref{eq:RemainderFunction}.
        \item If $U$ is uniform on $[-1,1]^{\Shell}$, then $\E_{U}(\norm{e(U)}_2^2)\le\rho_h^2\,s_*^2M$, where $\rho_h\coloneqq C_\rho(h+1)^{-(d-2)}$ with a constant $C_\rho$ depending only on $d$.
        \item $c_d(h+1)^{d-1/2}\le s_*M\le sM\le C_d(h+1)^{d-1/2}$.
    \end{enumerate}
\end{proposition}

\begin{proof}
    (a) By Fubini's theorem and \eqref{eq:Primitive},
    \begin{equation*}
        \inp{G_z}{W(u)}=\sum_i(Bz)_i\int_{\cube}\one_{\{P_i\le x\}}W(u;x)\,dx=\sum_i(Bz)_i\,(\Jst W(u))(P_i)=z\cdot V(u).
    \end{equation*}

    (b) For $f\in L^2(\cube)$, $\Jst\Plo f\in\Low$, so $\bigl((\Jst\Plo f)(P_i)\bigr)_i=\Lambda b$ for some $b\in\R^{\dim\Low}$, and $B^*\Lambda b=(\Lambda^*B)^*b=0$ by \cref{lem:BProperties}(a). Since $\Jst f=\Jst\Plo f+\Jst\Phig f$, this proves the second claim, and applying it to $f\in\Hlo$ (where $\Phig f=0$) together with the computation in (a) gives $\inp{G_z}{f}=0$, that is, $\Plo G_z=0$.

    (c) As in (a), $\inp{G_z}{q_R}=\sum_i(Bz)_i\,(\Jst q_R)(P_i)=v\sum_i(Bz)_iT_R(P_i)=v\,(T^*Bz)_R=\alpha\gamma_dv^2z_R$ by \cref{lem:BProperties}(b).

    (d) By (b), $V(u)=B^*\bigl((\Jst\Phig W(u))(P_i)\bigr)_i=B^*\bigl(\Flin(u)(P_i)\bigr)_i+e(u)$. By \eqref{eq:WlinFlin}, $\bigl(\Flin(u)(P_i)\bigr)_i=vT\diag(\lambda_R)u$, and $B^*\bigl(vT\diag(\lambda_R)u\bigr)=v(T^*B)^*\diag(\lambda_R)u=\alpha\gamma_dv^2\diag(\lambda_R)u$ by \cref{lem:BProperties}(b).

    (e) By \eqref{eq:BDefinition}, we have $e(u)=\alpha\gamma_dv\,\Gram^{-1}T^*E\Pi\,E\bigl(\Delta_u(P_i)\bigr)_i$. Since $\Pi^2=\Pi=\Pi^*$, we have $\norm{T^*E\Pi}_{\mathrm{op}}^2=\norm{T^*E\Pi ET}_{\mathrm{op}}=\norm{\Gram}_{\mathrm{op}}\le C_\star\gamma_dv$, and $\norm{\Gram^{-1}}_{\mathrm{op}}\le(c_\star\gamma_dv)^{-1}$ by \cref{lem:Gram}, so that
    \begin{equation*}
        \norm{e(u)}_2^2\le\alpha^2\gamma_d^2v^2\cdot\frac{C_\star\gamma_dv}{c_\star^2\gamma_d^2v^2}\sum_ia_i|\Delta_u(P_i)|^2=\frac{C_\star}{c_\star^2}\,\alpha^2\gamma_dv\sum_ia_i|\Delta_u(P_i)|^2 .
    \end{equation*}
    Taking expectations and using \eqref{eq:SamplingScalar} and \eqref{eq:ScalarIntegral},
    \begin{equation*}
        \E_{U}\paren[\big]{\norm{e(U)}_2^2}\le\frac{C_\star}{c_\star^2}\,\alpha^2\gamma_dv\sum_ia_iw(P_i)\le\frac{C_\star^2}{c_\star^2}\,\alpha^2\gamma_dv\int_{\cube}w\le C_d\,\alpha^2v^3\lambda^2(h+1)^{d-1}\theta^4 .
    \end{equation*}
    On the other hand, $s^2M=\alpha^2\gamma_d^2v^4\lambda^2K_h/v\ge c_d\alpha^2v^3\lambda^2(h+1)^{d-1}$ by \eqref{eq:ShapeCount}. Hence $\E_{U}(\norm{e(U)}_2^2)\le C_d\theta^4s^2M\le4C_d\theta^4s_*^2M=\rho_h^2s_*^2M$ with $C_\rho\coloneqq2C_d^{1/2}t_d^2$, where $C_d$ denotes the constant in the preceding inequality and we used $\theta^4=t_d^4(h+1)^{-2(d-2)}$; since $\alpha^2$ has cancelled, $C_\rho$ does not depend on $\omega_d$.

    (f) By \eqref{eq:BDefinition}, \eqref{eq:MDefinition} and \eqref{eq:SamplingSize}, $sM=\omega_d\sqrt{n/v}\,\gamma_dv\lambda K_h$ with $n\asymp M=K_h/v$, so $sM\asymp\sqrt{K_h}\,v^{-1}\cdot v\lambda K_h=\lambda K_h^{3/2}\asymp\lambda(h+1)^{3(d-1)/2}$ by \eqref{eq:ShapeCount}. Inserting $\lambda=\varepsilon t_d(h+1)^{-(d-2)/2}$ gives $sM\asymp(h+1)^{d-1/2}$, and $s_*\ge s/2$.
\end{proof}

\section{The entropy obstruction}
\label{sec:Entropy}

We now study the standard Gaussian vector $z\in\R^M$ conditioned on the event that all weighted box counts are small. For $t>0$ let
\begin{equation}\label{eq:SmallBall}
    K_t\coloneqq\bigl\{z\in\R^M:\norm{G_z}_\infty\le t\bigr\},
    \qquad
    \norm{G_z}_\infty\coloneqq\sup_{x\in\cube}|G_z(x)|,
\end{equation}
with $G_z$ as in \eqref{eq:IncidenceField}. With the enumeration of $\Shell$ fixed in \cref{sec:Witness}, $\R^M=\R^{\Shell}$ carries the coordinates $z_1,\dots,z_M$.

\begin{lemma}[Geometry of $K_t$]
    \label{lem:SmallBallSet}
    For every $t>0$, the set $K_t$ is closed, convex, symmetric and bounded, and it contains a neighborhood of the origin. In particular $0<\gm(K_t)$.
\end{lemma}

\begin{proof}
    For fixed $x$, the map $z\mapsto G_z(x)$ is linear, and $G_z(x)$ takes only finitely many values as $x$ ranges over the cube, namely the values $\sum_{i\in S}(Bz)_i$ over the finitely many incidence patterns $S=\{i:P_i\le x\}$. Hence $z\mapsto\norm{G_z}_\infty$ is a seminorm on $\R^M$, and $K_t$ is closed, convex and symmetric; being the unit ball of a seminorm on a finite-dimensional space, it contains a neighborhood of the origin. By \cref{prop:Response}(c), $\alpha\gamma_dv^2|z_R|=|\inp{G_z}{q_R}|\le\norm{G_z}_\infty\norm{q_R}_1=v\norm{G_z}_\infty$, so $K_t\subseteq[-t/(\alpha\gamma_dv),t/(\alpha\gamma_dv)]^M$.
\end{proof}

\subsection{Conditioned Gaussian vectors and Rosenblatt coordinates}
\label{sec:ConditionedGaussian}

Let $K\subset\R^M$ be a bounded Borel set with $p\coloneqq\gm(K)>0$, and let $\nu\coloneqq\gm(\,\cdot\mid K)$ be the conditioned Gaussian measure. It has the density
\begin{equation}\label{eq:ConditionedDensity}
    \varphi(z)=\frac{(2\pi)^{-M/2}\mathrm{e}^{-\norm{z}_2^2/2}}{p}\,\one_K(z),
\end{equation}
which is bounded and vanishes outside the bounded set $K$. Hence \cref{lem:ChainRule} applies to a random vector $Z$ with law $\nu$, and
\begin{equation}\label{eq:ConditionedEntropy}
    \ent(\nu)=-\int\varphi\log\varphi=\frac M2\log(2\pi)+\frac12\E_{\nu}\paren[\big]{\norm{Z}_2^2}+\log p\ \ge\ \frac M2\log(2\pi)+\log p .
\end{equation}
We only use the lower bound in \eqref{eq:ConditionedEntropy}: since $\ent(\gm)=\frac M2\log(2\pi\mathrm{e})$, the differential entropy of the conditioned vector is at most $\log(1/p)+M/2$ below that of $\gm$, that is, at most $\kappa+\frac12$ nats per coordinate when $p\ge\mathrm{e}^{-\kappa M}$. Let $U=U(Z)$ be the Rosenblatt transform of $Z$ under $\nu$, as in \cref{lem:Rosenblatt}: $U_j\coloneqq2F_j(Z_j\mid Z_{<j})-1$, where $F_j(\,\cdot\mid z_{<j})$ is the conditional distribution function of $Z_j$ given $Z_{<j}=z_{<j}$ under $\nu$.

\begin{proposition}[Correlation with the Rosenblatt coordinates]
    \label{prop:Correlation}
    Let $K$, $p$, $\nu$, $Z$ and $U$ be as above. Then $U$ is uniformly distributed on $[-1,1]^M$ under $\nu$, $\E_{\nu}(Z_jU_j)\ge0$ for every $j\le M$, and
    \begin{equation}\label{eq:CorrelationBound}
        \sum_{j=1}^M\E_{\nu}(Z_jU_j)\ \ge\ \frac{M}{4\mathrm{e}}\exp\Bigl(\frac{\ent(\nu)}{M}\Bigr)\ \ge\ \frac{\sqrt{2\pi}}{4\mathrm{e}}\,p^{1/M}\,M .
    \end{equation}
    In particular, if $p\ge\mathrm{e}^{-\kappa M}$ for some $\kappa>0$, then $\sum_j\E_{\nu}(Z_jU_j)\ge c_\kappa M$ with $c_\kappa\coloneqq\frac{\sqrt{2\pi}}{4\mathrm{e}}\mathrm{e}^{-\kappa}$.
\end{proposition}

\begin{proof}
    The uniform law of $U$ is \cref{lem:Rosenblatt}. We first prove a scalar inequality. Let $Y$ be a real random variable with a bounded density of bounded support, let $F$ be its distribution function, and let $Y'$ be an independent copy. Since $F$ is continuous, $\E(\operatorname{sgn}(Y-Y')\mid Y)=F(Y)-(1-F(Y))=2F(Y)-1$ almost surely, and exchanging the roles of $Y$ and $Y'$,
    \begin{equation}\label{eq:GiniIdentity}
        \E(Y\,(2F(Y)-1))=\E\paren[\big]{Y\operatorname{sgn}(Y-Y')}=\tfrac12\E\paren[\big]{(Y-Y')\operatorname{sgn}(Y-Y')}=\tfrac12\E\paren[\big]{|Y-Y'|}\ge0 .
    \end{equation}
    Conditioning on $Y$ and using Jensen's inequality, $\E(|Y-Y'|)\ge\E(|Y-\E(Y')|)=\E(|Y-\E Y|)$, and \cref{lem:MaxEntropy} gives $\E(|Y-\E Y|)\ge\mathrm{e}^{\ent(Y)}/(2\mathrm{e})$. Hence
    \begin{equation}\label{eq:ScalarCorrelation}
        \E(Y\,(2F(Y)-1))\ \ge\ \frac{\mathrm{e}^{\ent(Y)}}{4\mathrm{e}} .
    \end{equation}

    Now fix $j$. For $\nu_{<j}$-almost every $z_{<j}$, where $\nu_{<j}$ is the law of $Z_{<j}$ (namely for those $z_{<j}$ at which the marginal density of $Z_{<j}$ is positive), the conditional law of $Z_j$ given $Z_{<j}=z_{<j}$ has a bounded density of bounded support, namely $\varphi_{\le j}(z_{<j},\cdot)/\varphi_{<j}(z_{<j})$ in terms of the marginal densities of $\varphi$. Applying \eqref{eq:GiniIdentity} and \eqref{eq:ScalarCorrelation} to this conditional law gives
    \begin{equation*}
        \E_{\nu}(Z_jU_j\mid Z_{<j}=z_{<j})\ \ge\ \frac{1}{4\mathrm{e}}\exp\bigl(\ent(Z_j\mid Z_{<j}=z_{<j})\bigr)\ \ge\ 0 .
    \end{equation*}
    Taking the expectation over $z_{<j}$ and applying Jensen's inequality to the convex function $\exp$ yields $\E_{\nu}(Z_jU_j)\ge\frac{1}{4\mathrm{e}}\exp\bigl(\ent(Z_j\mid Z_{<j})\bigr)$, with the averaged conditional entropy of \cref{lem:ChainRule}. Summing over $j$, applying Jensen's inequality once more, and using the chain rule,
    \begin{equation*}
        \sum_{j=1}^M\E_{\nu}(Z_jU_j)\ \ge\ \frac{M}{4\mathrm{e}}\exp\Bigl(\frac1M\sum_{j=1}^M\ent(Z_j\mid Z_{<j})\Bigr)=\frac{M}{4\mathrm{e}}\exp\Bigl(\frac{\ent(\nu)}{M}\Bigr).
    \end{equation*}
    The second inequality in \eqref{eq:CorrelationBound} follows from \eqref{eq:ConditionedEntropy}.
\end{proof}

\subsection{The small-ball obstruction}
\label{sec:Obstruction}

We combine \cref{prop:Correlation} with the exact linear response of \cref{prop:Response}. Recall that $s_*=\min_Rs_R$ and that $\rho_h=C_\rho(h+1)^{-(d-2)}$ is the quantity in \cref{prop:Response}(e).

\begin{proposition}[Small-ball obstruction]
    \label{prop:Obstruction}
    Let $\kappa>0$ and $t>0$ be such that $\gm(K_t)\ge\mathrm{e}^{-\kappa M}$. Then
    \begin{equation*}
        t\ \ge\ (c_\kappa-\rho_h)\,s_*M,\qquad c_\kappa=\frac{\sqrt{2\pi}}{4\mathrm{e}}\mathrm{e}^{-\kappa} .
    \end{equation*}
\end{proposition}

\begin{proof}
    Let $\nu\coloneqq\gm(\,\cdot\mid K_t)$, let $Z$ have law $\nu$, and let $U=U(Z)$ be its Rosenblatt transform; by \cref{lem:SmallBallSet}, $K_t$ satisfies the assumptions of \cref{sec:ConditionedGaussian}. For every $z\in K_t$ and every $u\in[-1,1]^{\Shell}$, \cref{prop:Response}(a) and \cref{lem:Normalization} give
    \begin{equation*}
        z\cdot V(u)=\inp{G_z}{W(u)}\le\norm{G_z}_\infty\norm{W(u)}_1\le t .
    \end{equation*}
    Since $Z\in K_t$ almost surely and $U\in[-1,1]^{\Shell}$, it follows that $t\ge\E_{\nu}(Z\cdot V(U))$; the expectation exists because $Z$ is bounded on $K_t$ and $V$ is a polynomial on a compact set. By \cref{prop:Response}(d),
    \begin{equation*}
        \E_{\nu}(Z\cdot V(U))=\sum_{R\in\Shell}s_R\,\E_{\nu}(Z_RU_R)+\E_{\nu}(Z\cdot e(U)) .
    \end{equation*}
    Each term $\E_{\nu}(Z_RU_R)$ is nonnegative by \cref{prop:Correlation}, so the first sum is at least
    \begin{equation*}
        s_*\sum_{R\in\Shell}\E_{\nu}(Z_RU_R)\ \ge\ s_*c_\kappa M .
    \end{equation*}
    For the second term, the Cauchy--Schwarz inequality, \cref{lem:SecondMoment} and the fact that $U$ has exactly the uniform law on $[-1,1]^{\Shell}$ under $\nu$ give
    \begin{equation*}
        \bigl|\E_{\nu}(Z\cdot e(U))\bigr|\le\bigl(\E_{\nu}\paren[\big]{\norm{Z}_2^2}\bigr)^{1/2}\bigl(\E_{\nu}\paren[\big]{\norm{e(U)}_2^2}\bigr)^{1/2}\le\sqrt M\cdot\bigl(\E_{U}\paren[\big]{\norm{e(U)}_2^2}\bigr)^{1/2}\le\sqrt M\cdot\rho_hs_*\sqrt M ,
    \end{equation*}
    by \cref{prop:Response}(e). Hence $t\ge s_*c_\kappa M-\rho_hs_*M$. Note that no independence between $Z$ and $e(U(Z))$ is used, and that the Rosenblatt transform is never differentiated.
\end{proof}

\section{Proof of the main theorem}
\label{sec:Conclusion}

It remains to show that the hereditary discrepancy $D$ of the point set gives $\gm(K_{2D})\ge\mathrm{e}^{-\kappa M}$ with $\kappa$ depending only on $d$ (Step~1 of the overview), and to pass from the constructed point sets to all sufficiently large $N$.

\begin{lemma}[Gaussian small balls from hereditary discrepancy]
    \label{lem:GaussianCover}
    Let $A\in\R^{m\times n}$ satisfy $\herdisc(A)\le D$, and let $\eta$ be a centered Gaussian vector in $\R^n$ with $\Var(\eta_i)\le1$ for all $i\le n$. Then
    \begin{equation*}
        \Prob(\norm{A\eta}_\infty\le2D)\ge(p_0/2)^n,\qquad p_0\coloneqq\Prob(|g|\le1),
    \end{equation*}
    where $g$ is a standard Gaussian variable.
\end{lemma}

The argument is a Gaussian analogue of the volume argument behind the determinant lower bound of Lov\'asz, Spencer and Vesztergombi \cite[Lemma~2 and Theorem~2]{LSV86}.

\begin{proof}
    By \cref{thm:Sidak}, $\Prob(\eta\in[-1,1]^n)\ge\prod_i\Prob(|\eta_i|\le1)\ge p_0^n$, because $\Prob(|\eta_i|\le1)\ge p_0$ for a centered Gaussian variable of variance at most one. By \eqref{eq:RoundingCover}, the event $\{\eta\in[-1,1]^n\}$ is contained in the union over $\sigma\in\{-1,1\}^n$ of the events $\{A\eta\in A\sigma+K\}$ with $K\coloneqq[-2D,2D]^m$. The vector $A\eta$ is a centered Gaussian vector in $\R^m$, and $K$ is closed, convex and symmetric, so \cref{thm:Anderson} gives $\Prob(A\eta\in A\sigma+K)\le\Prob(A\eta\in K)$ for every $\sigma$. By the union bound over the $2^n$ vectors $\sigma$, $p_0^n\le\Prob(\eta\in[-1,1]^n)\le2^n\,\Prob(A\eta\in K)$.
\end{proof}

Let $h\in3\N$, let $P=\{P_1,\dots,P_n\}\subset(0,1)^d$ be the point set of \cref{sec:Sampling} at depth $h$, and let $A_P$ be the incidence matrix of the set system $\{P\cap A(x):x\in\cube\}$, as in \cref{sec:HereditaryDiscrepancy}. Put
\begin{equation}\label{eq:DhDefinition}
    \mathfrak{D}_h\coloneqq\herdisc(A_P).
\end{equation}

\begin{proposition}[Hereditary discrepancy of the construction]
    \label{prop:HerdiscLowerBound}
    There are constants $c_d>0$ and $h_d\in3\N$ such that the point set of \cref{sec:Sampling} at depth $h$ satisfies $\mathfrak{D}_h\ge c_d(h+1)^{d-1/2}$ for all $h\in3\N$ with $h\ge h_d$.
\end{proposition}

\begin{proof}
    Let $z$ be standard Gaussian in $\R^M$ and $\eta\coloneqq Bz$. By \cref{lem:BProperties}(c), $\Var(\eta_i)=(BB^*)_{ii}\le1$. The rows of $A_P$ are the vectors $(\one_{\{P_i\le x\}})_{i\le n}$, $x\in\cube$, so $\norm{A_P\eta}_\infty=\sup_x|G_z(x)|=\norm{G_z}_\infty$ and $\{\norm{A_P\eta}_\infty\le2\mathfrak{D}_h\}=\{z\in K_{2\mathfrak{D}_h}\}$. Note that $\mathfrak{D}_h\ge1$, so that $t=2\mathfrak{D}_h>0$ is admissible in \cref{prop:Obstruction}: the restriction of $A_P$ to a single column has discrepancy one, because the point $P_i$ lies in $A(x)$ for $x=(1,\dots,1)$. \Cref{lem:GaussianCover} with $D=\mathfrak{D}_h$ and \eqref{eq:SamplingSize} give
    \begin{equation*}
        \gm(K_{2\mathfrak{D}_h})\ge(p_0/2)^n=\mathrm{e}^{-a_0n}\ge\mathrm{e}^{-\kappa_dM},\qquad a_0\coloneqq\log(2/p_0),\quad \kappa_d\coloneqq a_0C_d,
    \end{equation*}
    where $C_d$ is the constant with $n\le C_dM$ in \eqref{eq:SamplingSize}. \Cref{prop:Obstruction} with $t=2\mathfrak{D}_h$ and $\kappa=\kappa_d$ yields $2\mathfrak{D}_h\ge(c_{\kappa_d}-\rho_h)s_*M$. Here $c_{\kappa_d}>0$ depends only on $d$, while $\rho_h=C_\rho(h+1)^{-(d-2)}$ decreases to zero as $h\to\infty$ because $d\ge3$; choose $h_d\in3\N$ such that $\rho_h\le c_{\kappa_d}/2$ for all $h\ge h_d$. For $h\ge h_d$ we obtain $2\mathfrak{D}_h\ge\frac12c_{\kappa_d}s_*M$, and \cref{prop:Response}(f) gives $\mathfrak{D}_h\ge\frac14c_{\kappa_d}s_*M\ge c_d(h+1)^{d-1/2}$.
\end{proof}

\begin{proof}[Proof of \cref{thm:Main}]
By \eqref{eq:SamplingSize}, \eqref{eq:MDefinition} and \eqref{eq:ShapeCount}, the construction at depth $h$ has $n\le C_d'(h+1)^{d-1}2^h$ points, with a constant $C_d'\ge1$. For sufficiently large $N$, let $h=h(N)$ be the largest element of $3\N$ with $C_d'(h+1)^{d-1}2^h\le N$. Then $h\le\log_2N$, and the maximality of $h$ gives $C_d'(h+4)^{d-1}2^{h+3}>N$, hence $h>\log_2N-\log_2(8C_d'(\log_2N+4)^{d-1})\ge\frac12\log_2N$ for all $N\ge N_d$, where $N_d$ is chosen so large that for every $N\ge N_d$ the number $h(N)$ is defined, $8C_d'(\log_2N+4)^{d-1}\le\sqrt N$ and $\log_2N\ge2h_d$; in particular $h(N)\ge h_d$.

For such $N$, the point set $P$ at depth $h(N)$ has $n\le N$ points, and \cref{prop:HerdiscLowerBound} together with \eqref{eq:HerdiscPoints} provides a subset $P'\subseteq P$ such that every coloring $\chi$ of $P'$ satisfies $\sup_{x\in\cube}|\sum_{p\in P'\cap A(x)}\chi(p)|\ge\mathfrak{D}_h\ge c_d(h+1)^{d-1/2}\ge c_d'(\log N)^{d-1/2}$ for a constant $c_d'>0$ depending only on $d$.

Let $P''$ be the union of $P'$ with $N-|P'|$ distinct points all of whose coordinates lie in $(1,2)$. For every $x\in\cube$ we have $P''\cap A(x)=P'\cap A(x)$, since the added points have no coordinate at most one, and every coloring of $P''$ restricts to a coloring of $P'$; as the supremum defining $\disc(P'',\Ad)$ runs over all $x\in\R^d$, in particular over $x\in\cube$, we get $\disc(P'',\Ad)\ge c_d'(\log N)^{d-1/2}$. Since $\disc(N,\Ad)$ is the supremum over all $N$-point sets in $\R^d$, this proves \cref{thm:Main} with the constant $c_d'$. Moreover, dilating $P''$ by $1/2$ gives an $N$-point set in $(0,1)^d$ satisfying the same lower bound.
\end{proof}

\begin{remark}[Order of the constants]
    \label{rem:Constants}
    For fixed $d$, the constants are chosen in the following order, none of them depending on $h$: $\varepsilon=\frac14$ and $t_d=t_d^*$ (\cref{thm:Remainder}); the sampling constants (\cref{lem:SamplingLemma}), which together with $t_d$ fix $C_\rho$ and hence the sequence $(\rho_h)_h$; $\omega_d$ (\cref{lem:BProperties}), which fixes $\alpha$, hence $s$ and $s_*$ for every $h$, and the constants in \cref{prop:Response}(f); the entropy constant $\kappa_d=a_0C_d$ and hence $c_{\kappa_d}$; finally $h_d$ and $c_d$ (\cref{prop:HerdiscLowerBound}), and then $N_d$ and $c_d'$ (\cref{thm:Main}). Although $\kappa_d$ depends on the sampling constant $C_d$ through $n\le C_dM$, there is no circularity: the sampling constants are not adapted to $\kappa_d$, because $\rho_h$ decays polynomially in $h$ for every fixed $d\ge3$, while $c_{\kappa_d}$ is a fixed positive number.
\end{remark}

\section{Extensions and further consequences}
\label{sec:Consequences}

This section proves \cref{cor:Polytopes,cor:BrownianSheet} and gives an efficient coloring algorithm attaining the upper bound in \cref{cor:Main}. The algorithm is independent of \cref{thm:Main}; that theorem establishes the optimality of its guarantee for $d\ge3$.

\subsection{Fixed convex polytopes}
\label{sec:Polytopes}

Recall that $\mathcal T_B=\{tB+x:t>0,\ x\in\R^d\}$ denotes the family of homothetic copies of a convex polytope $B\subset\R^d$. For a finite set $P\subset\R^d$ and a set $F\subseteq\R^d$ we call $P\cap F$ the \emph{trace} of $F$ on $P$.

\begin{lemma}[Finite-trace transfer]
    \label{lem:Transfer}
    Let $d\ge1$, let $B\subset\R^d$ be a convex polytope with nonempty interior and let $P\subset\R^d$ be finite. There are an invertible linear map $L$ of $\R^d$ and $t>0$ such that for every $x\in\R^d$ there is $y\in\R^d$ with $L(P\cap A(x))=L(P)\cap(tB+y)$.
\end{lemma}

\begin{proof}
    The idea is to find a linear map under which, for the finitely many relevant differences $x-p$, membership in the tangent cone of $B$ at a vertex coincides with $p\le x$; choosing $t$ large then makes the facets of $B$ away from this vertex irrelevant.

    Translating $B$ does not change $\mathcal T_B$, so we may assume that the origin is a vertex of $B$. Write
    \begin{equation*}
        B=\{y\in\R^d:\ell_i(y)\ge0\ \text{for }i\le m,\ a_r(y)\le b_r\ \text{for }r\le m'\},
    \end{equation*}
    where the linear functionals $\ell_1,\dots,\ell_m$ define the facets of $B$ through the origin, and the linear functionals $a_1,\dots,a_{m'}$ with constants $b_1,\dots,b_{m'}$ define the remaining facets; $b_r>0$ because the origin belongs to $B$ but lies on none of the latter facets. The tangent cone $C\coloneqq\{y:\ell_i(y)\ge0\ \text{for }i\le m\}$ has nonempty interior because it contains $B$, and it is pointed: since all $b_r>0$, a vector $y\ne0$ with $\ell_i(y)=0$ for all $i$ would give a segment $\{sy:|s|\le s_0\}$ in $B$ lying in all facets through the vertex $0$, whose intersection is $\{0\}$. Moreover, $tB=\{y\in C:a_r(y)\le tb_r\ \text{for }r\le m'\}$ for $t>0$.

    \emph{A basis adapted to $C$.} Choose faces $\{0\}=F_0\subset F_1\subset\dots\subset F_d=C$ of $C$ with $\dim F_j=j$, which exist because $C$ is pointed and every proper face of a polyhedral cone is contained in a face of the next dimension, and choose $v_j$ in the relative interior of $F_j$ for $j=1,\dots,d$. Then $v_1,\dots,v_d$ is a basis of $\R^d$. Indeed, the face $F_{j-1}=F_j\cap H$ is cut out of $F_j$ by a supporting linear hyperplane $H$ of $F_j$ (it passes through the origin since $F_j$ is a cone), so that $\operatorname{span}F_{j-1}\subseteq H$; the relative interior point $v_j$ of $F_j$ lies in $H$ only if $F_j\subseteq H$, which is excluded by $\dim F_j>\dim F_{j-1}$. Hence $v_j\notin\operatorname{span}\{v_1,\dots,v_{j-1}\}$. For $\varepsilon>0$ let $T_\varepsilon$ be the invertible linear map $T_\varepsilon z\coloneqq\sum_{j=1}^d\varepsilon^{j-1}z_jv_j$. Since $v_j\in C$, all $\ell_i(v_j)$ are nonnegative, and since $v_d$ is an interior point of $C$, all $\ell_i(v_d)$ are positive. Put $k(i)\coloneqq\min\{j:\ell_i(v_j)>0\}$. Then $\ell_i(v_j)=0$ for $j<k(i)$, so that for every fixed $z\in\R^d$
    \begin{equation}\label{eq:TransferLimit}
        \varepsilon^{1-k(i)}\ell_i(T_\varepsilon z)=\sum_{j\ge k(i)}\varepsilon^{j-k(i)}z_j\ell_i(v_j)\longrightarrow\ell_i(v_{k(i)})\,z_{k(i)},\qquad\varepsilon\downarrow0.
    \end{equation}
    Moreover, every $j\in\{1,\dots,d\}$ is of the form $k(i)$: the face $F_{j-1}$ is the intersection of the facets of $C$ containing it, and $F_j\not\subseteq F_{j-1}$, so some facet $C\cap\{\ell_i=0\}$ contains $F_{j-1}$ but not $F_j$; then $\ell_i(v_{j'})=0$ for $j'<j$, while $\ell_i\ge0$ on $F_j$ and $\ell_i\not\equiv0$ on $F_j$ force $\ell_i(v_j)>0$ at the relative interior point $v_j$, that is, $k(i)=j$.

    \emph{Traces.} Let $X\subset\R^d$ be a finite set containing, for every trace $P\cap A(x)$ with $x\in\R^d$, one point $x$ with this trace, chosen such that $x_j\ne p_j$ for all $p\in P$ and $j\le d$; this is possible because increasing a coordinate of $x$ slightly does not change the trace. All coordinates of the finitely many differences $x-p$ with $x\in X$ and $p\in P$ are nonzero, so by \eqref{eq:TransferLimit} there is $\varepsilon>0$ such that $\ell_i(T_\varepsilon(x-p))$ has the sign of $(x-p)_{k(i)}$ for all $i$, $x$ and $p$. As the indices $k(i)$ exhaust $\{1,\dots,d\}$, we obtain
    \begin{equation*}
        T_\varepsilon(x-p)\in C\iff(x-p)_j>0\ \text{for all }j\le d\iff p\le x,\qquad p\in P,\ x\in X,
    \end{equation*}
    where the last equivalence uses $x_j\ne p_j$. Fix this $\varepsilon$, put $L\coloneqq-T_\varepsilon$, and choose $t>0$ with $a_r(L(p-x))\le tb_r$ for all $p\in P$, $x\in X$ and $r\le m'$. Then, for $p\in P$ and $x\in X$,
    \begin{equation*}
        L(p)\in tB+L(x)\iff L(p-x)\in C\iff T_\varepsilon(x-p)\in C\iff p\le x,
    \end{equation*}
    where the first equivalence uses the choice of $t$; that is, $L(P)\cap(tB+L(x))=L(P\cap A(x))$. Since every trace $P\cap A(x)$ with $x\in\R^d$ equals $P\cap A(x')$ for some $x'\in X$, the claim follows with $y\coloneqq L(x')$.
\end{proof}

\begin{proof}[Proof of \cref{cor:Polytopes}]
    Let $P\subset\R^d$ be any $N$-point set and let $L$ and $t$ be as in \cref{lem:Transfer}. The set $Q\coloneqq L(P)$ has $N$ points, the colorings of $Q$ are the maps $\chi\circ L^{-1}$ with $\chi$ a coloring of $P$, and by \cref{lem:Transfer}, $L$ maps every set $P\cap A(x)$ onto a set $Q\cap(tB+y)$ with $tB+y\in\mathcal T_B$. Hence $\disc(Q,\mathcal T_B)\ge\disc(P,\Ad)$, and taking the supremum over $P$ gives $\disc(N,\mathcal T_B)\ge\disc(N,\Ad)\ge c_d(\log N)^{d-1/2}$ for $N\ge N_d$ by \cref{thm:Main}. The upper bound $\disc(N,\mathcal T_B)=O_{d,B}((\log N)^{d-1/2})$ is \cite[Theorem~2]{Nikolov17}.
\end{proof}

For generic polytopes, \cite[Theorem~3]{Nikolov17} gives the lower bound $\Omega_{d,B}((\log N)^{d-1})$; Nikolov conjectured that the genericity assumption there is unnecessary and that the lower-bound constant can be chosen independently of $B$ \cite[p.~1093]{Nikolov17}. Both conjectures already follow from \cref{lem:Transfer}: the proof of \cref{cor:Polytopes} gives $\disc(N,\mathcal T_B)\ge\disc(N,\Ad)$ in every dimension, and $\disc(N,\Ad)=\Omega_d((\log N)^{d-1})$ for $d\ge2$ \cite{MatousekNikolov15,MNT20}. The bound for all sufficiently large $N$ follows from the grid estimate in \cite[Section~6]{MNT20} by taking a subset realizing hereditary discrepancy and padding, as in the proof of \cref{thm:Main}. For $d\ge3$, \cref{cor:Polytopes} improves the exponent to $d-\tfrac12$, with the constants $c_d,N_d$ of \cref{thm:Main}.

\subsection{A matching offline coloring guarantee}
\label{sec:Algorithm}

We use the following recent result on prefix discrepancy. For operation counts, we count arithmetic operations and comparisons on real numbers, evaluations of square roots and exponentials, and sampling from $\mathcal N(0,1)$ or $\operatorname{Uniform}[0,1)$ at unit cost.

\begin{literaturetheorem}[{Aden-Ali \cite[Theorem~1]{AdenAli26}}]
    \label{thm:AdenAli}
    Let $v_1,\dots,v_T\in\R^m$ be a fixed sequence with $\norm{v_i}_2\le1$ for all $i\le T$; in particular, the vectors do not depend on the random choices of the algorithm below. There is a randomized algorithm which receives $v_1,\dots,v_T$ one at a time and chooses a sign $\sigma_i\in\{-1,1\}$ upon receiving $v_i$, before seeing $v_{i+1},\dots,v_T$, with the following property: for every $k\le T$, the prefix sum $S_k\coloneqq\sum_{i\le k}\sigma_iv_i$ can be written as $S_k=G_{k,1}+G_{k,2}+G_{k,3}$ with random vectors $G_{k,1},G_{k,2},G_{k,3}$ each having the standard Gaussian distribution on $\R^m$ (they need not be independent). In the model above, the algorithm runs in time $O(mT)$.
\end{literaturetheorem}

By H\"older's inequality, $\E(\exp(\inp{S_k}y))\le\prod_{j=1}^3\bigl(\E(\exp(3\inp{G_{k,j}}y))\bigr)^{1/3}=\exp(9\norm{y}_2^2/2)$ for every $y\in\R^m$, so every prefix sum is subgaussian with variance proxy $9$, as noted in \cite[Section~1]{AdenAli26}, and the Chernoff bound gives, for $w\in\R^m\setminus\{0\}$ and $u>0$,
\begin{equation}\label{eq:PrefixTail}
    \Prob(|\inp{S_k}w|\ge u)\le2\exp\Bigl(-\frac{u^2}{18\norm{w}_2^2}\Bigr).
\end{equation}

We follow the reduction in Nikolov's proof of \cite[Theorem~1]{Nikolov17}, using an explicit dyadic factorization and replacing Banaszczyk's signed-series theorem \cite{Banaszczyk12} (\cite[Lemma~9]{Nikolov17}) with \cref{thm:AdenAli}.

\begin{proposition}[Efficient offline coloring for anchored boxes]
    \label{prop:Algorithm}
    Let $d\ge2$. There is a randomized offline algorithm (it reads the whole point set before choosing the first sign) which, given an $N$-point set $P\subset\R^d$ with $N\ge2$ and $\delta\in(0,1)$, outputs a coloring $\chi\colon P\to\{-1,1\}$ such that, with probability at least $1-\delta$,
    \begin{equation}\label{eq:AlgorithmBound}
        \sup_{x\in\R^d}\Bigl|\sum_{p\in P\cap A(x)}\chi(p)\Bigr|\le3K\sqrt{2\log\frac{2(N+1)^d}{\delta}},\qquad K\coloneqq(\lceil\log_2N\rceil+1)^{d-1},
    \end{equation}
    and on the same event the discrepancy of $\chi$ with respect to all axis-parallel boxes is at most $2^d$ times this bound. The algorithm performs $O_d(N^2(\log N)^{d-1})$ operations in the model of \cref{thm:AdenAli}.
\end{proposition}

\begin{proof}
    \emph{Reduction to ranks and prefixes.} Put $H\coloneqq\lceil\log_2N\rceil$, so that $N\le2^H$, and fix an order of the points of $P$. For $j\le d-1$ let $r_j(p)\in\{1,\dots,N\}$ be the rank of $p$ among the points of $P$ with respect to the $j$-th coordinate, ties broken by the fixed order, and let $r(p)\coloneqq(r_1(p),\dots,r_{d-1}(p))$. For every threshold $x_j$, the points with $p_j\le x_j$ are exactly the points with $r_j(p)\le\rho_j$, where $\rho_j\in\{0,\dots,N\}$ is their number, because tied points are included together. Enumerate the points as $P_1,\dots,P_N$ in nondecreasing order of the last coordinate, ties again broken by the fixed order; then the points $p$ with $p_d\le x_d$ are $P_1,\dots,P_k$ for some $k\in\{0,\dots,N\}$. Consequently every trace $P\cap A(x)$ has the form
    \begin{equation}\label{eq:TraceForm}
        \{P_i:i\le k,\ r(P_i)\le\rho\}\qquad\text{with }k\in\{0,\dots,N\}\text{ and }\rho\in\{0,\dots,N\}^{d-1},
    \end{equation}
    and there are at most $(N+1)^d$ distinct traces.

    \emph{Dyadic products.} For $0\le l\le H$ the dyadic integer intervals $\{(q-1)2^l+1,\dots,q2^l\}$, $1\le q\le2^{H-l}$, partition $\{1,\dots,2^H\}$, so every rank vector $r(p)$ lies in exactly $K=(H+1)^{d-1}$ products of dyadic intervals, one for each choice of the levels $l_1,\dots,l_{d-1}$. Let $\calQ$ be the set of dyadic products containing at least one rank vector $r(P_i)$; then $m\coloneqq|\calQ|\le NK$. For $i\le N$ let $v_i\in\R^{\calQ}$ be the indicator vector of the products in $\calQ$ containing $r(P_i)$, divided by $\sqrt K$, so that $\norm{v_i}_2=1$. The algorithm applies \cref{thm:AdenAli} to $v_1,\dots,v_N$ in this order and outputs $\chi(P_i)\coloneqq\sigma_i$.

    \emph{Discrepancy.} Fix a trace \eqref{eq:TraceForm}. Each initial segment $\{1,\dots,\rho_j\}$ of $\{1,\dots,2^H\}$ is the disjoint union of at most $H+1$ dyadic intervals, one for each binary digit of $\rho_j$, so $\{s\in\{1,\dots,2^H\}^{d-1}:s\le\rho\}$ is the disjoint union of at most $K$ dyadic products. Let $w\in\{0,1\}^{\calQ}$ be the indicator vector of those among them that belong to $\calQ$; then $\norm{w}_2^2\le K$. The traces and the vectors $w$ depend only on $P$, not on the random signs. For every point $P_i$ of the trace, $r(P_i)$ lies in exactly one of these products, and for no other $i\le k$ does $r(P_i)$ lie in any of them, so
    \begin{equation*}
        \sum_{i\le k,\ r(P_i)\le\rho}\sigma_i=\sqrt K\,\inp{S_k}w,\qquad S_k=\sum_{i\le k}\sigma_iv_i .
    \end{equation*}
    If $k=0$ or $w=0$, the signed count vanishes. Otherwise \eqref{eq:PrefixTail} gives $\Prob(\sqrt K|\inp{S_k}w|\ge u)\le2\exp(-u^2/(18K^2))$ for $u>0$. With $u\coloneqq3K\sqrt{2\log(2(N+1)^d/\delta)}$ this probability is at most $\delta(N+1)^{-d}$, and a union bound over the at most $(N+1)^d$ traces gives \eqref{eq:AlgorithmBound} with probability at least $1-\delta$. For an axis-parallel box in $\Bd$, the indicator function of its trace on $P$ is a signed combination of the indicator functions of $2^d$ traces $P\cap A(x)$ by inclusion--exclusion, as in \eqref{eq:AnchoredVsBoxes}, which gives the factor $2^d$ on the same event.

    \emph{Operations.} Sorting and ranking take $O(dN\log N)$ comparisons. The $NK$ incidences between points and dyadic products can be enumerated and indexed in $O_d(NK\log(NK))$ operations. Constructing the vectors $v_1,\dots,v_N$ and applying \cref{thm:AdenAli} takes $O(mN)=O(N^2K)=O_d(N^2(\log N)^{d-1})$ operations, since $m\le NK$, and this dominates. Since the ranks and the order of the points depend on the whole input, the algorithm is offline.
\end{proof}

For $\delta=N^{-a}$ with fixed $a>0$, the bound \eqref{eq:AlgorithmBound} is $O_{d,a}((\log N)^{d-1/2})$, since $K\asymp(\log N)^{d-1}$ and $\log(2(N+1)^d/\delta)=O_{d,a}(\log N)$. For $d\ge3$, \cref{thm:Main} shows that no coloring algorithm can guarantee a smaller order of worst-case discrepancy, so the discrepancy guarantee of \cref{prop:Algorithm} is optimal up to the constant. The operation count refers to the model of \cref{thm:AdenAli}; finite-precision (bit) complexity is not addressed.

\subsection{Small deviations of the Brownian sheet}
\label{sec:BrownianSheet}

The standard Brownian sheet $B_d$ on $\cube$ is the centered Gaussian process with continuous paths and covariance
\[
    \E(B_d(x)B_d(y))=\prod_{j=1}^d\min\{x_j,y_j\}.
\]

\begin{proof}[Proof of \cref{cor:BrownianSheet}]
The upper bound in \cref{cor:BrownianSheet} is due to
Dunker, Linde, K\"uhn and Lifshits \cite{DunkerLindeKuhnLifshits99}, for $d=2$ already to Bass \cite{Bass88},
and the lower bound for $d=2$ is due to Talagrand \cite{Talagrand94}.
We prove the lower bound for $d\ge3$.
Fix $h\in3\N$ and retain the notation of \cref{sec:Witness,sec:Remainder}.
Put
\[
    f_R\coloneqq\frac{T_R}{\sqrt{\gamma_dv}},\qquad
    (Jf)(x)\coloneqq\int_{[0,x]}f(p)\,dp,\qquad
    X_z\coloneqq\sum_{R\in\Shell}z_RJf_R.
\]
By \cref{lem:TOrthogonality} and \cref{lem:SpaceInclusions}(c), the $f_R$ are orthonormal
and orthogonal to $\Low=\Jst\Hlo$.
The operator $J$ maps $L^2(\cube)$ into $C(\cube)$, and Fubini gives
$\inp{Jf}{w}=\inp{f}{\Jst w}$ for $f,w\in L^2(\cube)$.
For $t>0$, the set $K_t^{(h)}\coloneqq\{z:\norm{X_z}_\infty\le t\}$ is
closed, convex and symmetric, and contains a neighborhood of zero.
It is bounded because
$\inp{X_z}{q_R}=\sqrt{\gamma_d}\,v^{3/2}z_R$ and $\norm{q_R}_1=v$ give
$|z_R|\le t/\sqrt{\gamma_dv}$ on $K_t^{(h)}$.

Let $\Delta_u=\Jst\Phig W(u)-\Flin(u)$ and define
$a_R\coloneqq\sqrt{\gamma_d}\,v^{3/2}\lambda_R$,
$a_*\coloneqq\min_Ra_R$, and
$r(u)\coloneqq(\inp{f_R}{\Delta_u})_{R\in\Shell}$.
Orthogonality to $\Low$ and \eqref{eq:WlinFlin} give
\begin{equation}\label{eq:BrownianPairing}
    \inp{X_z}{W(u)}=\sum_Rz_R\bigl(a_Ru_R+r_R(u)\bigr).
\end{equation}
For uniform $U$ on $[-1,1]^M$, Bessel's inequality, \cref{thm:Remainder}
and \eqref{eq:FlinNorm} imply
\begin{equation}\label{eq:BrownianRemainder}
    \E\paren[\big]{\norm{r(U)}_2^2}
    \le\E\paren[\big]{\norm{\Delta_U}_2^2}
    \le\frac{C_d'\theta^4}{3}\sum_Ra_R^2
    \le\frac{4C_d'}{3}\theta^4a_*^2M,
\end{equation}
where the last step uses $a_R\le2a_*$ from \cref{lem:LinearPart}.

Suppose $p\coloneqq\gm(K_t^{(h)})\ge\mathrm e^{-M}$, let
$Z\sim\gm(\,\cdot\mid K_t^{(h)})$, and let $U=U(Z)$ be its Rosenblatt
transform. By \cref{prop:Correlation}, $U$ is uniform,
$\E(Z_RU_R)\ge0$, and $\sum_R\E(Z_RU_R)\ge c_1M$, where
$c_1=\sqrt{2\pi}/(4\mathrm e^2)$.
Using \cref{lem:Normalization,lem:SecondMoment},
\eqref{eq:BrownianPairing}, \eqref{eq:BrownianRemainder} and
Cauchy--Schwarz, with all expectations under the conditioned law, yields
\[
    t\ge\E(\inp{X_Z}{W(U)})
    \ge c_1a_*M-\sqrt{M}\bigl(\E\paren[\big]{\norm{r(U)}_2^2}\bigr)^{1/2}
    \ge\left(c_1-2\sqrt{C_d'/3}\,\theta^2\right)a_*M.
\]
Since $\theta^2\to0$ and
$a_*M\asymp\lambda K_h\sqrt v\asymp 2^{-h/2}(h+1)^{d/2}$,
there is $b_d>0$ such that, for all sufficiently large $h\in3\N$ and
$g\sim\gm$,
\begin{equation}\label{eq:BrownianShellProbability}
    \Prob(\norm{X_g}_\infty\le b_d\delta_h)<\mathrm e^{-M},
    \qquad \delta_h\coloneqq2^{-h/2}(h+1)^{d/2}.
\end{equation}

To compare with $B_d$, realize it as
$B_d(x)=\mathcal W(\one_{[0,x]})$, where $\mathcal W$ is a centered Gaussian
linear family on $L^2(\cube)$ with covariance
$\E(\mathcal W(f)\mathcal W(w))=\inp{f}{w}$.
With $g_R\coloneqq\mathcal W(f_R)$ and $(Jf_R)(x)=\inp{f_R}{\one_{[0,x]}}$, we get
$B_d=X_g+Y_h$ with $Y_h(x)\coloneqq\mathcal W\bigl(\one_{[0,x]}-\sum_R(Jf_R)(x)f_R\bigr)$.
Since the $f_R$ are orthonormal, $g\sim\gm$, and $Y_h$ is independent of $g$,
because every argument of $\mathcal W$ in $Y_h$ is orthogonal to all $f_R$.
Both terms have continuous paths.
Apply \cref{thm:Anderson} on increasing finite subsets of a countable dense
set, condition on $Y_h$, and pass to the limit to obtain
\[
    \Prob(\norm{B_d}_\infty\le t)\le\Prob(\norm{X_g}_\infty\le t).
\]
Finally, $\delta_h\downarrow0$ eventually along $h\in3\N$, with consecutive
ratios bounded away from zero, and
\[
    M\asymp2^h(h+1)^{d-1}
    =\delta_h^{-2}(h+1)^{2d-1}
    \asymp\delta_h^{-2}\bigl(\log(1/\delta_h)\bigr)^{2d-1}.
\]
For every sufficiently small $\delta>0$, choose $h\in3\N$ with
$b_d\delta_{h+3}<\delta\le b_d\delta_h$.
Then $\delta\asymp\delta_h$, and monotonicity together with
\eqref{eq:BrownianShellProbability} gives the required bound.
\end{proof}

\section*{Acknowledgments}
The author thanks Aleksandar Nikolov for suggesting the Brownian-sheet application. The author used Anthropic's Claude and OpenAI's Codex to assist with literature research, proof development and checking, Lean formalization, and manuscript preparation. The author retains full responsibility for the content of the manuscript, including all mathematical claims and references.

\bibliographystyle{alpha}
\bibliography{References}
\end{document}